\newif\ifCONDITION
\CONDITIONtrue

\documentclass[preprint,12pt]{elsarticle}

\usepackage{setspace} 
\usepackage{graphicx,amsmath,latexsym,amssymb}
\usepackage{ascmac}
\usepackage{bm}
\usepackage{lscape}
\usepackage{mathrsfs}
\usepackage{subfigure}
\usepackage[usenames]{color}

 \usepackage{lineno}

\journal{Computational Physics}

\bmdefine{\bX}{X}
\bmdefine{\bx}{x}
\bmdefine{\btau}{\tau}
\bmdefine{\bn}{n}
\bmdefine{\bt}{t}
\bmdefine{\bb}{b}
\bmdefine{\bv}{v}
\bmdefine{\br}{r}
\bmdefine{\bS}{S}
\bmdefine{\bV}{V}
\bmdefine{\bFV}{FV}
\bmdefine{\be}{e}
\bmdefine{\bR}{R}
\bmdefine{\bE}{E}
\bmdefine{\bA}{A}
\bmdefine{\ba}{a}
\bmdefine{\bB}{B}
\bmdefine{\bzeta}{\zeta}
\bmdefine{\bU}{U}
\bmdefine{\bu}{u}
\bmdefine{\mT}{\mathcal{T}}
\bmdefine{\mG}{\mathcal{G}}

\newcommand{\del}{\partial}
\newcommand{\wid}{0.3}
\newcommand{\datahead}{100}
\newcommand{\timestep}{125}

\begin{document}

\begin{frontmatter}


\title{
A note on implementation of boundary variation diminishing algorithm to high-order local polynomial-based schemes
}

\author[label1]{Yoshiaki Abe}
\ead{abe.y.ax@m.titech.ac.jp}
\author[label1]{Ziyao Sun}
\author[label1]{Feng Xiao}

\address[label1]{Department of Mechanical Engineering, Tokyo Institute of Technology\\
4259 Nagatsuta Midori-ku, Yokohama, 226-8502, Japan.}
\setstretch{1.}
\begin{abstract}
\indent 
 A novel approach for selecting appropriate reconstructions is implemented to the hyperbolic conservation laws in the high-order local polynomial-based framework, e.g., the discontinuous Galerkin (DG) and flux reconstruction (FR) schemes.
The high-order polynomial approximation generally fails to correctly capture a strong discontinuity inside a cell due to the Runge phenomenon, which is replaced by more stable approximation on the basis of a troubled-cell indicator such as that used with the total variation bounded (TVB) limiter.
This paper examines the applicability of a new algorithm, so-called boundary variation diminishing (BVD) reconstruction, to the weighted essentially non-oscillatory (WENO) methodology in the FR framework including the nodal type DG method. 
The BVD reconstruction adaptively chooses a proper approximation for the solution function so as to minimize the jump between values at the left- and right-side of cell boundaries.

Using the BVD algorithm, several numerical tests are conducted for selecting an appropriate function between the original high-order polynomial and the WENO reconstruction approximating either smooth or discontinuous profiles.
The selected functions based on the BVD algorithm are not always the same as those by the conventional TVB limiter, which indicates that the present BVD algorithm offers a radically new criteria for selecting a reconstruction function without any ad hoc TVB parameter.
Subsequently, the computation of a linear advection equation is examined for the BVD and TVB criteria with the WENO methodology in the FR framework.
The results of the BVD algorithm are comparable to those using the conventional TVB limiter in terms of oscillation suppression and numerical dissipation, which would be also possible in the system equations such as Euler equations.
Furthermore, since the present BVD algorithm does not need any ad hoc constant such as the TVB parameter, it could be more reliable than the conventional TVB limiter that is prevailing in the existing DG and FR approaches for shocks and other discontinuities.
Note that the present work is limited to third or lower order polynomials, leaving the implementations for higher-order schemes a future work.
\end{abstract}
%
\vspace{-1cm}
\begin{keyword}
boundary variation diminishing; BVD; discontinuity; shock capturing; TVB limiter; WENO; DG; FR; troubled-cell indicator.
\end{keyword}
\end{frontmatter}

%
\renewcommand\thefootnote{\roman{footnote})}
\setstretch{1.2}
\section{Introduction}
Discretization of hyperbolic conservation laws with high-order spatial accuracy has been one of the most demanded topics specifically for accurate and efficient computation of flows.
Recently, several high-order local polynomial-based schemes, e.g., discontinuous Galerkin (DG)\cite{Cockburn1989}\cite{Cockburn1998}, spectral difference (SD), spectrral volume (SV), and flux reconstruction (FR/CPR) schemes\cite{Huynh2007}, have been demonstrated to show their applicability to the discretization of compressibile Euler / Navier-Stokes equations around complex geometries\cite{Wang2012e}.
The important characteristic of these schemes is multiple degrees of freedom in one cell, where high-order polynomials are defined in each cell.
These high-order polynomials are discontinuous between adjacent cells since they are individually constructed in each cell so that some flux evaluation method such as an approximate Riemann solver is required.
Therefore, these schemes can be regarded as high-resolution Godunov type methods.

The use of high-order polynomials is, of course, suitable for the representation of solution profile with spatially high wave-number such as turbulent flows.
In such a case, for example, the DG scheme satisfies a local cell entropy inequality which ensures a nonlinear stablity for the polynomial approximation with an arbitrary degree~\cite{Jiang1996}.
\footnote{Additional stabilization techniques have been also proposed on the basis of preservation of high-order moments such as kinetic energy~\cite{Jameson2008} and entropy~\cite{Ismail2009}, which have been originally developped in the finite-difference framework~\cite{Arakawa1966}\cite{Morinishi2010}\cite{Pirozzoli2011}.}.
However, since these high-order unstructured schemes are inherently linear schemes in the sense of Godunov \cite{Godunov1959}, some nonlinear limiting procedure is required to capture discontinuities avoiding the Runge phenomenon which by nature appears in high-order polynomial approximation~\cite{Shu1988}\cite{Shu1989}.
The typical ideas for these limiting procedure in the high-order unstructured schemes are as follows:
1) the troubled cells are detected by an indicator based on the unlimited solution, where the solution polynomial in the corresponding cell needs to be limited;
2) the high-order polynomial approximation (unlimited solution) is replaced by nonoscillatory solution using nonlinear recosntruction techniques.
Excessive efforts have been devoted to propose the appropriate troubled-cell indicator and nonlinear reconstruction methods for the unstructured schemes.
The mile-stone work can be seen in the combination of a total variation bounded limiter (so-called TVB limiter) and weighted essentially non-oscillatory (WENO) reconstruction\cite{Shu1988}\cite{Jiang1996}\cite{Liu1994}\cite{Shu2009} by Qiu et al.\cite{Qiu2005} and Zhong et al.\cite{Zhong2013} for the DG framework, which has been subsequently extended to the FR (CPR) schemes\cite{Du2015}.
Although the combination of TVB limiter and WENO reconstruction is successfully applied to practical flow computations as well as benchmark tests, it is well known that an ad hoc constant, so-called the TVB parameter $M$, remains undetermined in the TVB limiter so that several possible solutions exist for the shock related problems (see details in Qiu et al.\cite{Qiu2005} and Zhong et al.\cite{Zhong2013}).
Note that another novel approach using subcell limiting and reconstructions has been recently proposed in the DG framework by Dumbser et al.~\cite{Dumbser2014}, which would be more suitable for these local polynomial-based schemes.

This paper represents an implementation of a novel algorithm for appropriate reconstructions, so-called the boundary variation diminishing (BVD) algorithm~\cite{Sun2016}, to the hyperbolic conservation laws in the FR framework\cite{Huynh2007}. 
The BVD algorithm adaptively chooses appropriate reconstructions so as to minimize the jump between the values at the both side of cell boundaries, which does not require any ad hoc constants such as the conventional TVB parameter.
We focus on the FR schemes in this note since they can recover another local polynomial-based schemes, e.g., the nodal DG scheme, by adopting particular correction functions for linear problems. 

The rest of this note is organized as follows.
Section \ref{sec_NA} describes a numerical algorithm including the FR scheme and limiting procedures based on the existing TVB limiter and the new BVD reconstruction.
In Sec.\ref{sec_NT}, several numerical tests are conducted:
the selection of appropriate reconstructions for given profiles;
the computation of a linear advection equation for a discontinuous initial condition, where the applicability of the new BVD algorithm is examined.
\section{Numerical algorithm}\label{sec_NA}
\subsection{Discretization of governing equations in the FR framework}
Let us consider a one-dimensional hyperbolic conservation law as
\begin{align}
\frac{\del u}{\del t} + \frac{\del f}{\del x}=0,
\end{align}
where $u=u(x,t)$ is the solution function, and $f=f(u)$ is the flux function.
The physical domain $x\in [x_{\text{min}},x_{\text{max}}]$ is spatially subdivided into $N$ cells at each time step.
The $i$th cell is expressed as $I_i:=\{x | x_{\text{min};i}\leq x \leq  x_{\text{max};i}\}$, which is mapped to a standard cell $E_s:=\{\xi | -1\leq \xi \leq 1\}$.
Hereinafter, the polynomial approximation of an arbitrary function $\psi(\xi)$ is denoted as 
\begin{align}
I[\psi](\xi):=\sum_{i=0}^K\psi(\xi_i)\phi_i(\xi),\quad \phi_{i}(\xi)=\prod_{j=0,j\neq i}^{K}\frac{\xi-\xi_i}{\xi_j-\xi_i},
\end{align}
where $\phi_i(\xi)$ is the $K$th-order Lagrange polynomial.
The solution at the $i$th cell is approximated by the $K$th-order polynomial of $\xi$, i.e., $I[u_i](\xi)$, which is updated by
%
%
\begin{align}
\frac{\del I[u_i]}{\del t}+
\frac{\del I[f_{i}]}{\del \xi}+\left[f^{com}_{i-1/2}- I[f_{i}](-1)\right]\frac{d g_-}{d \xi}+\left[f^{com}_{i+1/2}-I[f_{i}](+1)\right]\frac{d g_+}{d \xi}=0,\label{Gov}
\end{align}
Here, $i\pm 1/2$ indicate the cell boundary of the $i$th cell, which are the same as the notation of $\xi=\pm 1$;
$f^{com}_{i\pm 1/2}$ are the so-called common flux shared between adjacent cells, which are usually determined by an approximate Riemann solver;
$g_{L}$ ($g_R$) is the so-called correction function, which is the $(K+1)$th-order polynomial and returns $1$ and $0$ at the left and right (right and left) boundaries of the $i$th cell, respectively.
The second and third terms in Eq.\eqref{Gov} modify the original discontinuous polynomial of $I[f_i](\xi)$ to be continuous between neighbouring cells with respect to the left and right boundaries, respectively.
Therefore, the conservation of $u$, i.e., the integral of $u$ over the entire computational domain, is ensured.

Although the polynomial approximation $I[u_i](\xi)$ of degree $K$ provides a high-resolution reconstruction for smoothed flows, it is not suitable for representing the profiles including strong discontinuities inside the cell since high-order polynomials near the discontinuitiy oscillate due to the Runge phenomenon.
Therefore, if the high-order polynomial approximation $I[u_i](\xi)$ is not appropriate in the $i$th cell, the cell should be marked as a troubled one, where the high-order polynomial is replaced by more robust and nonoscillatory one.
The troubled-cell indicator is conventionally constructed using the TVB limiter as will be described in Sec.\ref{sec_TVB}.
Our new idea contributes to such an appropriate reconstruction without existing TVB limiters (see Sec.\ref{sec_BVD}).

\subsection{Existing TVB limiter}\label{sec_TVB}
One of the most prevailing troubled-cell indicator is that used with the TVB limiter\cite{Qiu2005}\cite{Zhong2013}\cite{Shu2009}, which is overviewed in this subsection.
Denote the volume-integrated average (VIA) of the solution $u$ in the $i$th cell by
\begin{align}
\overline{u}_i:=\frac{1}{\Delta \xi}\int_{E_s}I[u_i]d\xi=\frac{1}{2}\int_{-1}^{1}I[u_i]d\xi.
\end{align}
The difference of the VIA and solution at cell boundaries are defined as
\begin{align}
\tilde{u}^{+}_i&:=u_{i+1/2}-\overline{u}_i=I[u_i](+1)-\overline{u}_i,\\
\tilde{u}^{-}_i&:=\overline{u}_i-u_{i-1/2}=\overline{u}_i-I[u_i](-1).
\end{align}
Then, $\tilde{u}^{\pm}_i$ are modified by the TVB modified minmod function\cite{Zhong2013} $\tilde{m}$ as
\begin{align}
\tilde{u}^{\pm;\text{mod}}_i&:=\tilde{m}(\tilde{u}^{\pm}_i,\Delta_+\overline{u}_i,\Delta_-\overline{u}_i),\\
\Delta_+\overline{u}_i:=&\overline{u}_{i+1}-\overline{u}_{i},\quad
\Delta_-\overline{u}_i:=\overline{u}_{i}  -\overline{u}_{i-1}.
\end{align}
The function $\tilde{m}$ is generally defined by
\begin{align}
&\tilde{m}(a_1,\ldots,a_l)=
\begin{cases}
a_1\quad &\text{if } |a_1|\leq Mh^2,\\
m(a_1,\ldots,a_l)&\text{otherwise},
\end{cases}\\
&m(a_1,\ldots,a_l)=
\begin{cases}
s \min |a_l|\quad &\text{if } s=\text{sign}(a_1)=\cdots=\text{sign}(a_l),\\
0&\text{otherwise}.
\end{cases}
\end{align}
where $h$ is the maximum cell interval ($h=\text{max}\Delta x_i$).
The ad hoc TVB parameter $M$ has to be fixed adequately depending on the problems.
If the minmod function $\tilde{m}$ returns the first argument, i.e., $\tilde{u}^{\pm;\text{mod}}_i=\tilde{u}^{\pm}_i$, the cell is not a troubled cell and original high-order polynomial $I[u_i](\xi)$ can be adopted.
Otherwise, the cell is marked as the troubled cell and another approximation such as the WENO reconstruction should replace the original one.
In this way, the TVB limiter can be employed as a troubled-cell indicator leaving the TVB parameter $M$ undetermined.

In summary, the solution $u_i(\xi)$ is switched as follows:
\begin{align}
u_i(\xi)=
\begin{cases}
u^{<1>}_i(\xi)=I[u_i](\xi)&\text{for nontroubled cells},\\
u^{<2>}_i(\xi)&\text{for troubled cells},
\end{cases}\label{eq_selection}
\end{align}
where $u^{<1>}_i(\xi)$ is a generally high-order polynomial one (original approximation), which is a $K$th-order polynomial in this note;
$u^{<2>}_i(\xi)$ is of the low order, which is more stable and less oscillatory one.
In this paper, we adopt WENO~\cite{Qiu2005} and simple WENO~\cite{Zhong2013} reconstructions for the second candidate as $u^{<2>}_i:=u^{WENO}_i\quad\text{or}\quad u^{SWENO}_i$ \footnote{Another possible candidate is the THINC (Tangent of Hyperbola for INterface Capturing) function which has been adopted for the BVD framework in a finite-volume fashion\cite{Sun2016}.}.
The detail of the WENO and simple WENO reconstructions is not presented here as it is not a main focus of the present paper.

\subsection{New reconstruction based on the BVD algorithm}\label{sec_BVD}
In Eq.\eqref{eq_selection}, the first or second candidate $u^{<p>}_i(\xi)\quad (p=1,2)$ approximating the solution function is selected depending on the troubled-cell markers, i.e., $u^{<2>}_i(\xi)$ being selected in the troubled cells and $u^{<1>}_i(\xi)$ otherwise.
In this section, we describe a guideline for the BVD reconstruction that adaptively chooses appropriate reconstructions so as to minimize the jump between the values at the left and right side of cell boundaries.
The algotrithm is as follows:
\begin{itemize}
\item[i)]
Prepare two reconstructions $u^{<1>}_i(\xi)$ and $u^{<2>}_i(\xi)$ in each cell.
$u^{<1>}_i(\xi)$ is the original high-order polynomial $I[u_i](\xi)$;
$u^{<2>}_i(\xi)$ is more stable (generally, lower-order) reconstruction, that is set to be the WENO or simple WENO reconstructions in this paper.
\item[ii)]
At the cell boundary of $i+1/2$ between $I_i$ and $I_{i+1}$, find $u_i^{<p>}(\xi)$ and $u_{i+1}^{<p>}(\xi)$ with $p$ and $q$ being either $1$ or $2$, so that the boundary variation (BV):
\begin{align}
BV(u)_{i+1/2}:=|u_{i}^{<p>}(+1)-u_{i+1}^{<q>}(-1)|\label{eq_BV1}
\end{align}
is minimized.
Note that $u_{i}^{<p>}(+1)$ and $u_{i+1}^{<q>}(-1)$ indicate the right boundary value of $u_{i}^{<p>}(\xi)$ at $i$th cell and left boundary value of $u_{i+1}^{<q>}(\xi)$ at $(i+1)$th cell, respectively.
\item[iii)]
When the selected reconstructions are inconsistent between left and right boundaries of $i$th cell, that is, $u^{<p>}_i(\xi)$ found to minimize
\begin{align}
BV(u)_{i-1/2}:=|u_{i-1}^{<q>}(+1)-u_{i}^{<p>}(-1)|\label{eq_BV2}
\end{align}
is different from that found to minimize the BV value of \eqref{eq_BV1}, the following condition is adopted to uniquely determine the reconstruction function:
\begin{align}
u_{i}^{<p>}(\xi)=
\begin{cases}
u_i^{<1>}(\xi),\quad &\text{if } (\overline{u}_i-\overline{u}_{i+1})(\overline{u}_{i-1}-\overline{u}_{i})<0,\\
u_i^{<2>}(\xi),\quad &\text{otherwise.}
\end{cases}
\end{align}
\end{itemize}
The third step can be replaced by the criterion based on the boundary values as was introduced in the finite-volume framework\cite{Sun2016}.

\section{Numerical tests}\label{sec_NT}
First, the selection of appropriate reconstructions based on Sec.\ref{sec_TVB} and \ref{sec_BVD} are examined.
Then, the applicability of the new criterion for a linear advection equation is shown.
The numerical schemes examined in this section are listed in Table \ref{Tab_scheme}.
\begin{table}[htbp]\label{Tab_scheme}
{\footnotesize
\setstretch{1.5}
\caption{Numerical schemes}\vspace{0.1cm}
\begin{tabular}{ccc}
\hline\hline
 Troubled-cell indicator & First candidate: $u^{<1>}_i$ & Second (stable) candidate: $u^{<2>}_i$\\\hline
TVB & Normal $K$th-order polynomial: $I[u_i]$ & WENO: $u^{WENO}_i$\\
TVB & Normal $K$th-order polynomial: $I[u_i]$ & Simple WENO: $u^{SWENO}_i$\\
BVD & Normal $K$th-order polynomial: $I[u_i]$ & WENO: $u^{WENO}_i$\\
BVD & Normal $K$th-order polynomial: $I[u_i]$ & Simple WENO: $u^{SWENO}_i$\\\hline\hline
\end{tabular}
\setstretch{1.2}
}
\end{table}
The solution points are located at the Gauss points of $K+1$;
the correction functions $g_{\pm}$ are set to be the Radau polynomial, i.e., $g_{DG}$ \cite{Huynh2007};
the flux evaluation is the first-order upwind scheme at each cell boundaries;
the three-stage TVD Runge-Kutta scheme is used for time marching with a CFL number of $0.01$.

\subsection{Polynomial selection by the TVB and BVD criteria}
Here, we examine the selection of polynomial candidates $u^{<1>}_i$ and $u^{<2>}_i$ approximating the given analytical solution on the basis of the TVB and BVD criteria.
The computational domain is set to be $0\leq x\leq 10$;
the number of total cells are $N=3$.
The analytical solutions are given as:
\begin{itemize}
\item Case 1: smoothed flow
\begin{align}
u&=A\sin(2\pi\omega x),\label{eq_Case1}\\
A&=0.20,\quad \omega=0.25,
\end{align}
\item Case 2: discontinuous flow
\begin{align}
u=
\begin{cases}
1.0&\quad (x\leq 5.0),\label{eq_Case2}\\
0.0&\quad ( 5.0<x).
\end{cases}
\end{align}
\end{itemize}
In the smooth profiles of Eq.\eqref{eq_Case2}, the number of solution points per wavelength is equal to or more than three if $K\geq 2$.
\subsubsection{Results of the case 1: smooth flow}
Figure \ref{fig_Case1_K2} shows the polynomial approximation of degree $K=2$ for the smooth flow given by Eq.\eqref{eq_Case1}.
The grey lines show the analytical solution;
black lines show the high-order polynomial approximation $u^{<1>}_i=I[u_i]$;
blue lines show the stable function $u^{<2>}_i=u^{WENO/SWENO}_i$;
red lines show the selected function based on TVB or BVD criteria.
The selection of functions in the first and third cell are fixed to the high-order one $u_i=u^{<1>}_i$ for ease of understandings.
In Figs.\ref{fig_Case1_K2}(a)--(c), the WENO reconstruction $u^{<2>}_i$ in the second cell (blue-colored lines) is almost linear and quite different from either the exact profile (grey-colored lines) or the high-order polynomial approximation $u^{<1>}_i$ (black-colored lines).
The important features are summarized as:
\begin{itemize}
\item The BVD criteria in Fig.\ref{fig_Case1_K2}(a) selects the high-order polynomial $u^{<1>}_i$ so as to minimize the jump between neighbouring cells, which would provide the better approximation for such a smooth flow.
\item The TVB limiter with $M=0$ in Fig.\ref{fig_Case1_K2}(b), on the other hand, selects the stable function $u^{<2>}_i$ since the $M=0$ results in possibly the most stable criteria (TVD).
When $M$ is increased as $M=200$ in Fig.\ref{fig_Case1_K2}(c), the high-order polynomial $u^{<1>}_i$ is selected, which would be better than the case with $M=0$.
\end{itemize}
In this way, the classical TVB limiter contains uncertain properties of $M$ regarding the selection of functions, which does not always show better results even for the smooth flows.
Note that an optimum $M$ values can be estimated from the smoothness of the initial solution\cite{Shu1988}, which is not generally suitable for a nonlinear problem and not considered herein.
The simple WENO reconstruction $u^{<2>}_i$ in Figs.\ref{fig_Case1_K2}(d)--(f) corresponds to the high-order polynomial $u^{<1>}_i$.

\subsubsection{Results of the case 2: discontinuous flow}
Next, the discontinuous cases (Figs.\ref{fig_Case2_K2}(a)--(c)) shows that the WENO reconstruction $u^{<2>}_i$ in the second cell is a monotonically dicreasing profile. 
The important features are summarized as:
\begin{itemize}
\item
The TVB limiter with smaller $M$ (Fig.\ref{fig_Case2_K2}(b)) selects the stable WENO reconstruction $u^{<2>}_i$ while $M=200$ case (Fig.\ref{fig_Case2_K2}(c)) provides the original high-order polynomial $u^{<1>}_i$ although it contains a little overshoot and undershoot near the discontinuity.
\item
The BVD criteria in Fig.\ref{fig_Case2_K2}(a) selects more stable one $u^{<2>}_i$, where the condition Eq.\eqref{eq_BV2} is enacted.
On the other hand, the BVD criteria selects $u^{<1>}_i$ for the simple WENO reconstruction (Fig.\ref{fig_Case2_K2}(d)), which would result in higher-resolution solution although it may lead an unstable and oscillatory solution.
\end{itemize}
Therefore, it would be difficult to promise that the BVD criteria always selects more stable function than the TVB one only from the present observation.
However, the BVD criteria realizes a completely paramter-free selection, of which property would be better than that of the classical TVB limiter.
Note that practically, the BVD selection affects the polynomial selection in both of the left- and right-side of each cell boundaries.
The more complicated cases are discussed in Appendix A where the trend is almost the same as these basic cases.
\renewcommand{\datahead}{bv_data_cfl0p1_lhs3_rhs1_cor1_sp0_init}
\renewcommand{\timestep}{00000001}
\renewcommand{\wid}{0.3275} 
\begin{figure}[!htbp]
    \centering
    \setcounter{subfigure}{0}
    \subfigure[][{\scriptsize BVD/WENO}]{\includegraphics[width=\wid\textwidth,clip]{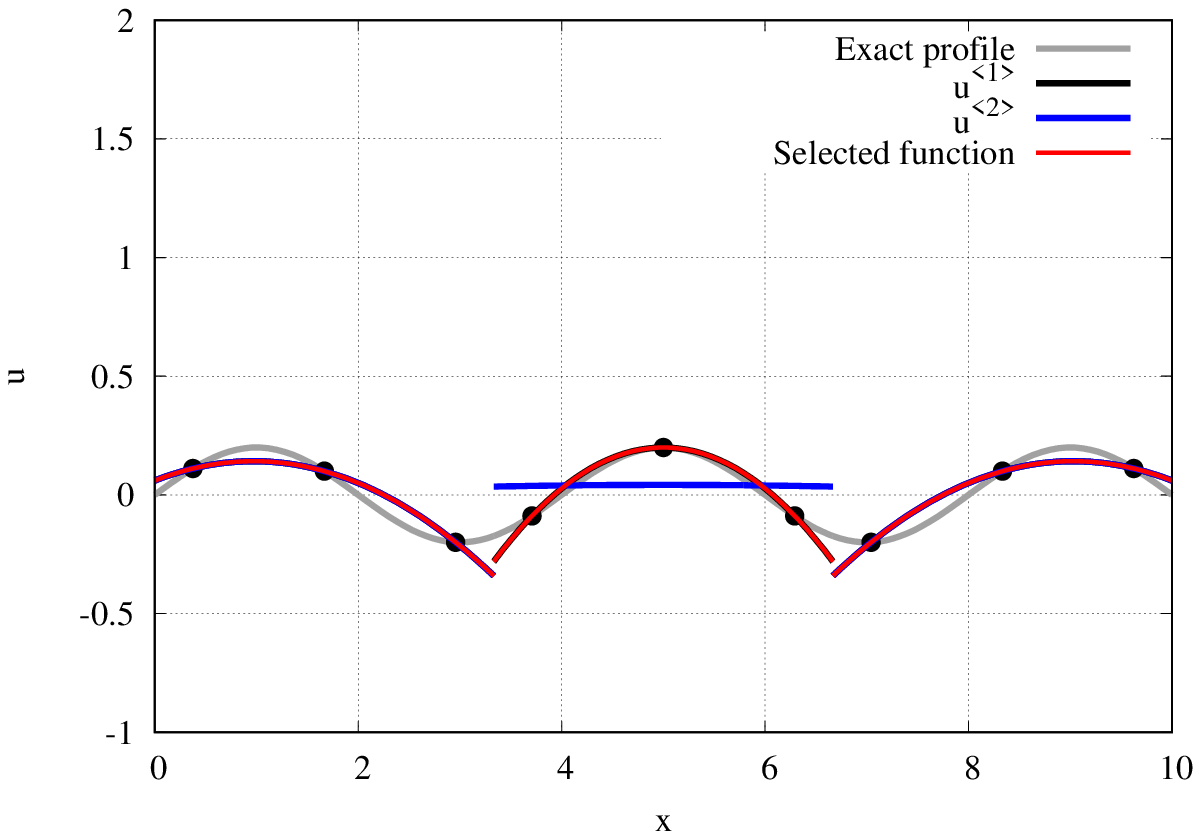}}
    \subfigure[][{\scriptsize TVB/WENO $(M=0.0)$}]{\includegraphics[width=\wid\textwidth,clip]{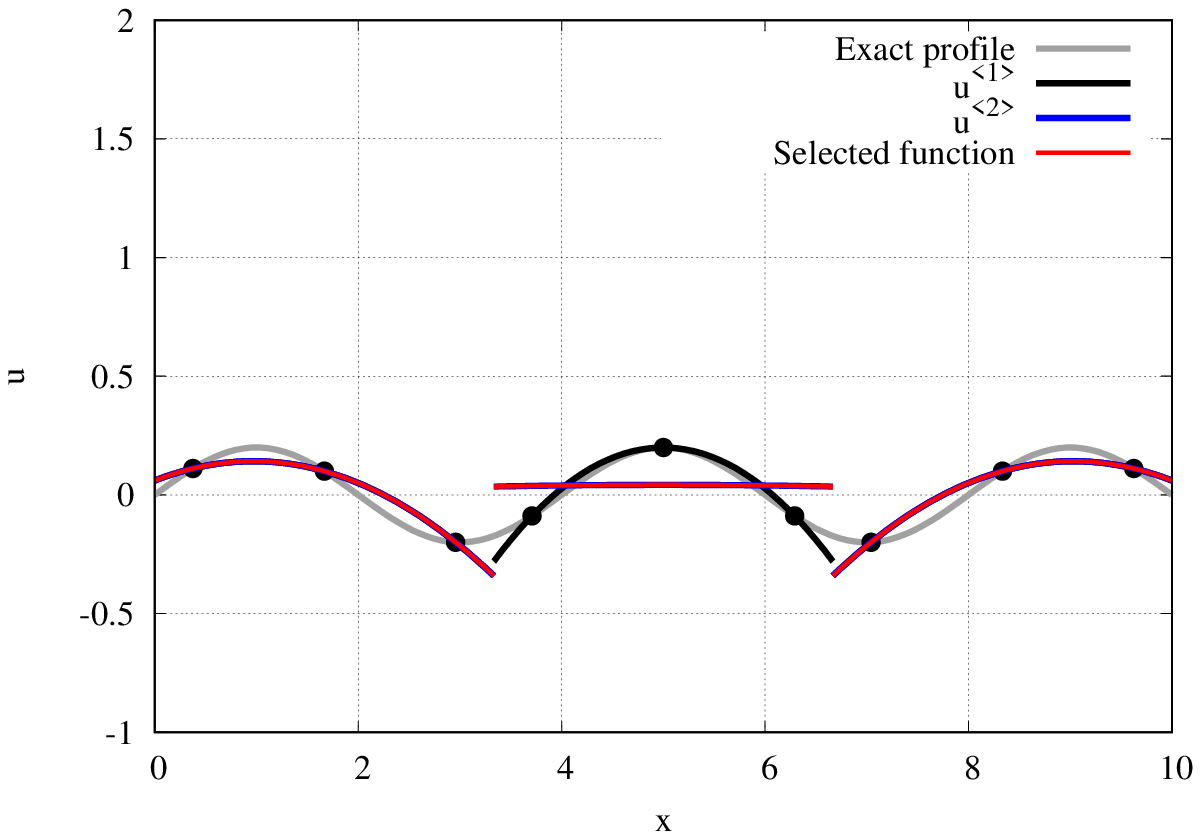}}
    \subfigure[][{\scriptsize TVB/WENO $(M=200)$}]{\includegraphics[width=\wid\textwidth,clip]{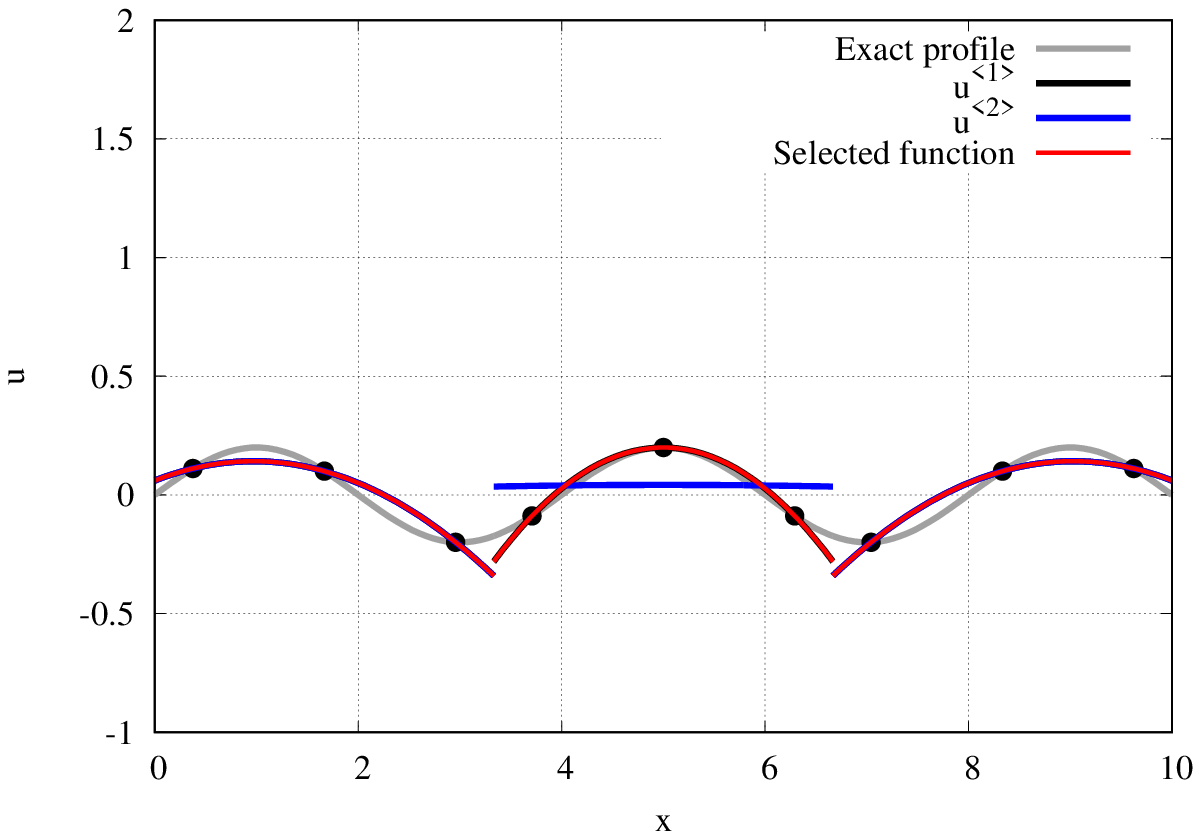}}\\
    \subfigure[][{\scriptsize BVD/SWENO}]{\includegraphics[width=\wid\textwidth,clip]{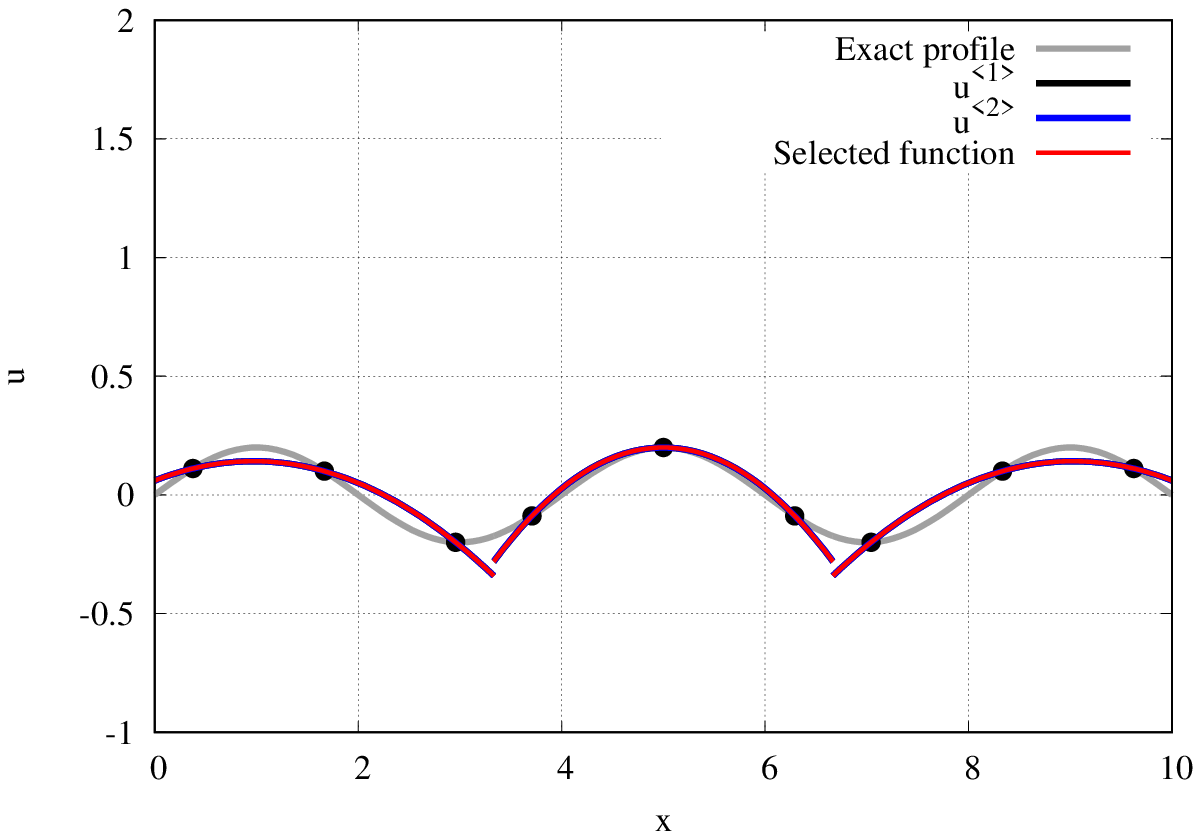}}
    \subfigure[][{\scriptsize TVB/SWENO $(M=0.0)$}]{\includegraphics[width=\wid\textwidth,clip]{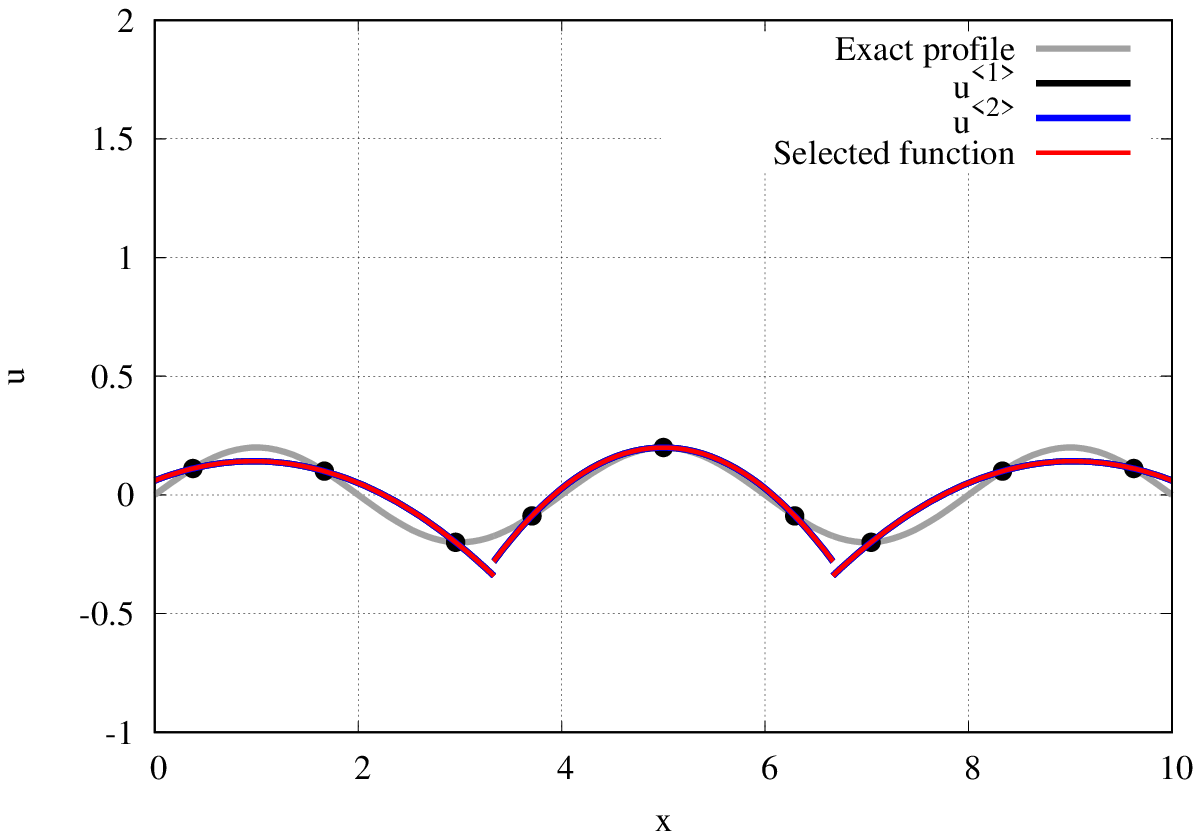}}
    \subfigure[][{\scriptsize TVB/SWENO $(M=200)$}]{\includegraphics[width=\wid\textwidth,clip]{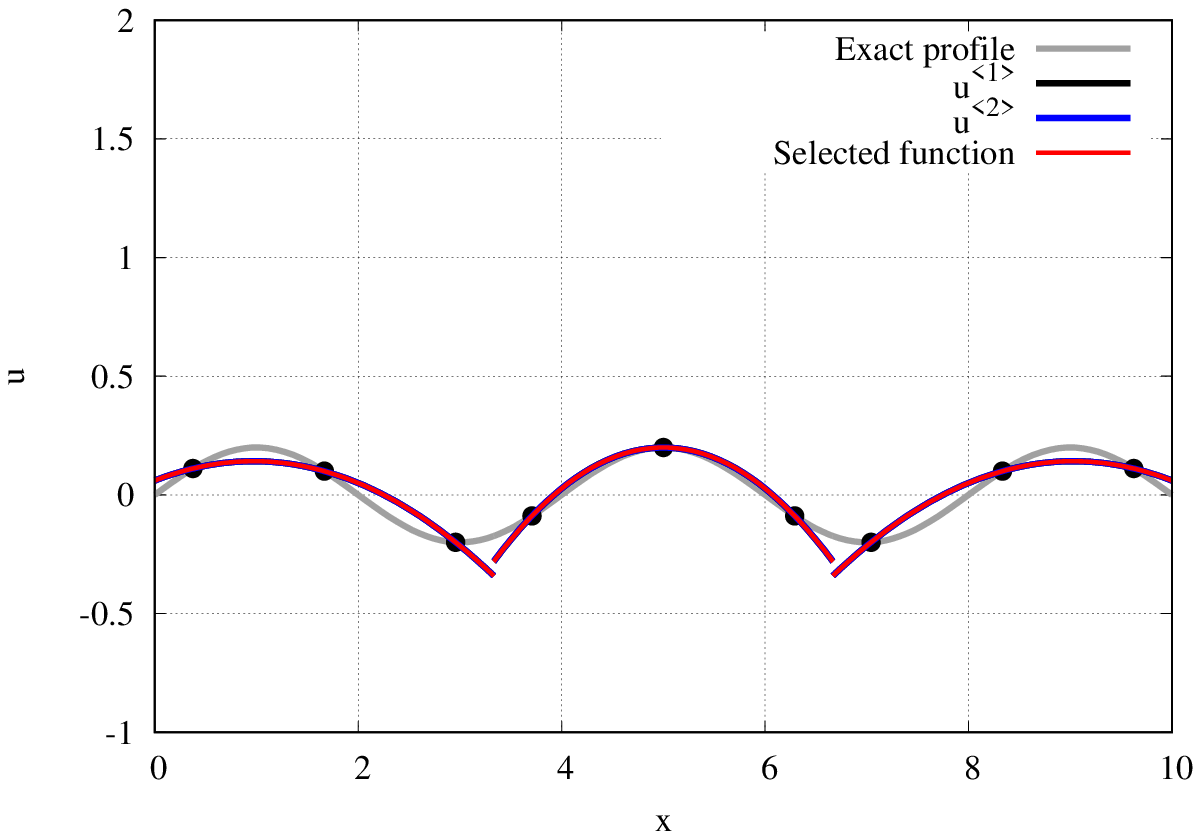}}
\caption{Polynomial selection for the Case 1 solution (smooth flow given by Eq.\eqref{eq_Case1}); polynomial degree $K=2$.
The grey lines show the analytical solution;
black lines show the high-order polynomial approximation $u^{<1>}_i$;
blue lines show the stable function $u^{<2>}_i$;
red lines show the selected function based on TVB or BVD criteria.
The computational domain is divided into three cells, where the functions in the first and third cell are fixed to $u^{<1>}_i$
 }\label{fig_Case1_K2}
\end{figure}
\begin{figure}[!htbp]
    \centering
    \setcounter{subfigure}{0}
    \subfigure[][{\scriptsize BVD/WENO}]{\includegraphics[width=\wid\textwidth,clip]{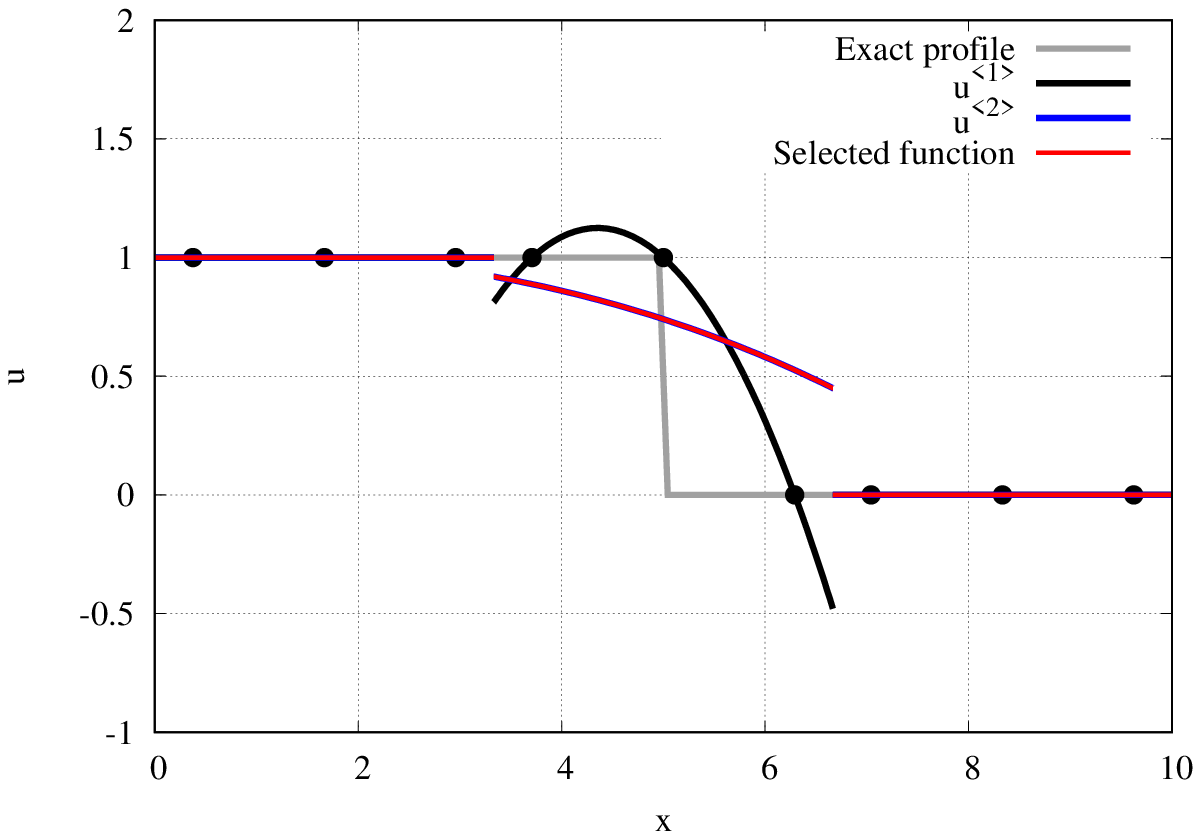}}
    \subfigure[][{\scriptsize TVB/WENO $(M=0.0)$}]{\includegraphics[width=\wid\textwidth,clip]{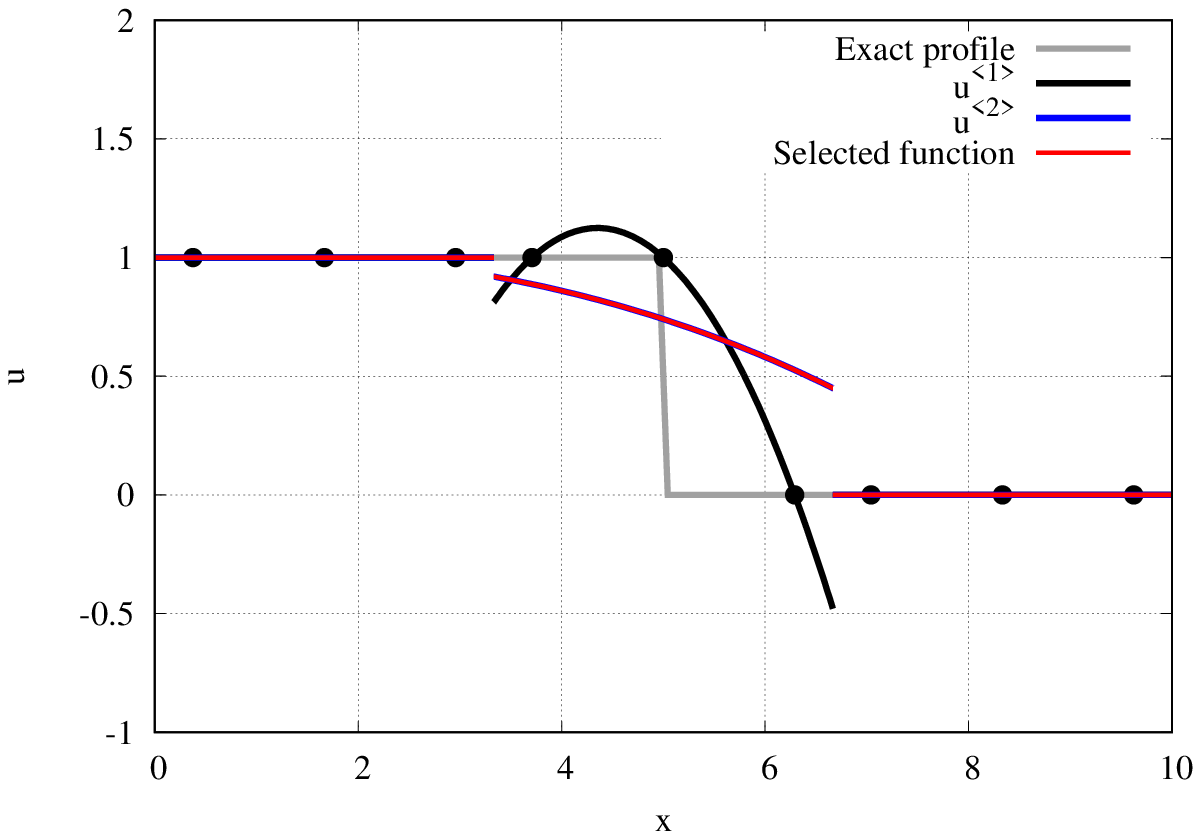}}
    \subfigure[][{\scriptsize TVB/WENO $(M=200)$}]{\includegraphics[width=\wid\textwidth,clip]{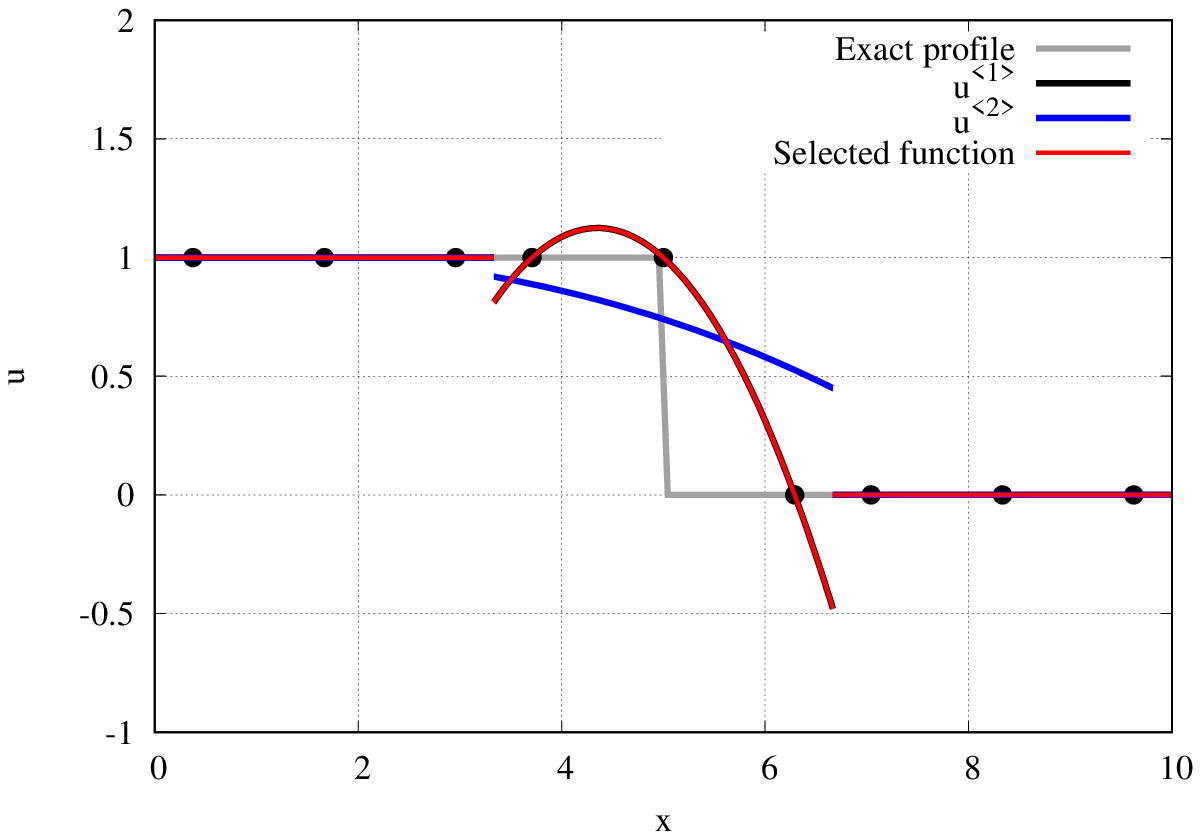}}\\
    \subfigure[][{\scriptsize BVD/SWENO}]{\includegraphics[width=\wid\textwidth,clip]{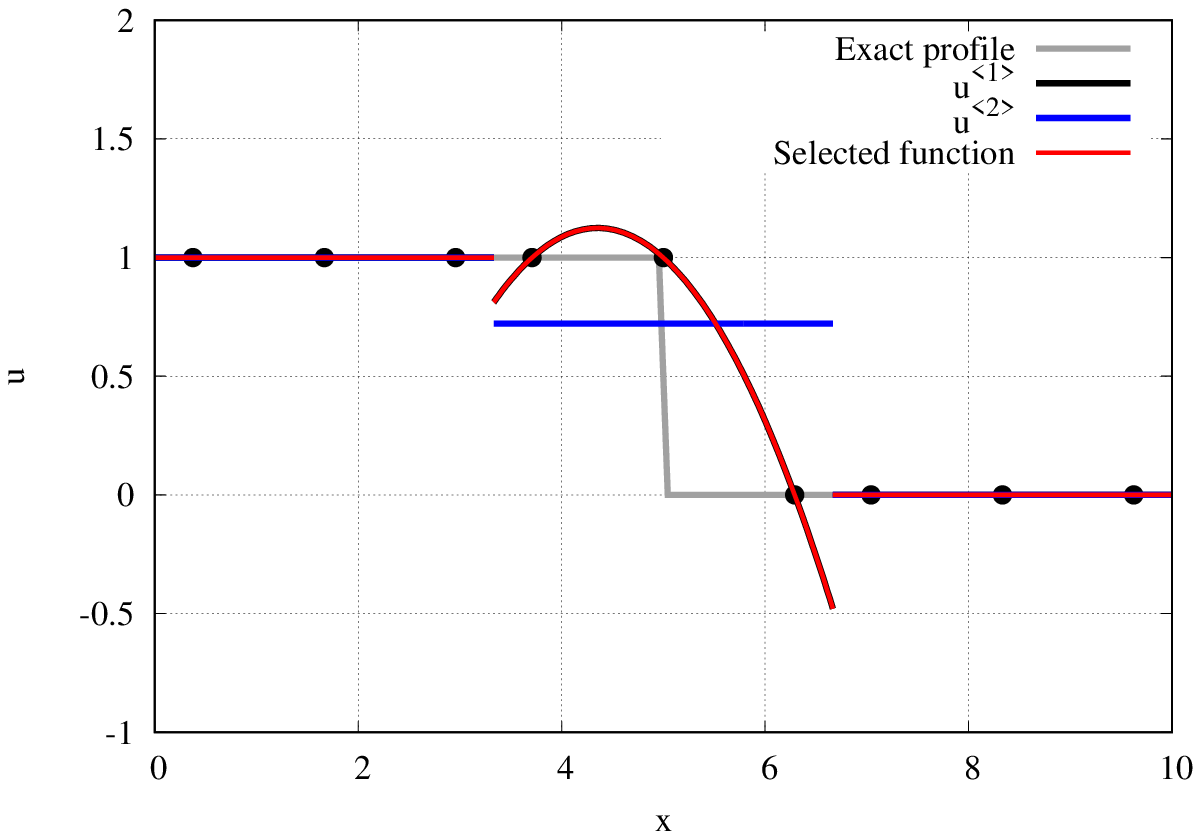}}
    \subfigure[][{\scriptsize TVB/SWENO $(M=0.0)$}]{\includegraphics[width=\wid\textwidth,clip]{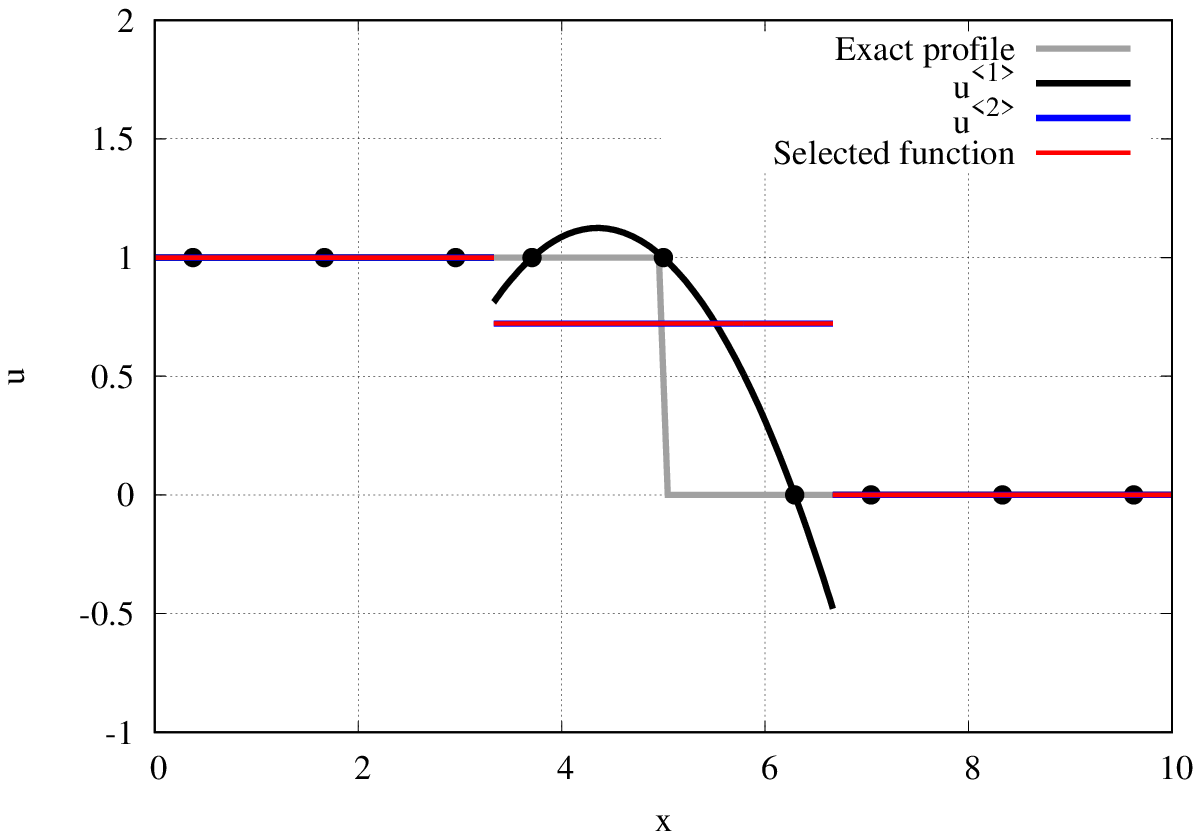}}
    \subfigure[][{\scriptsize TVB/SWENO $(M=200)$}]{\includegraphics[width=\wid\textwidth,clip]{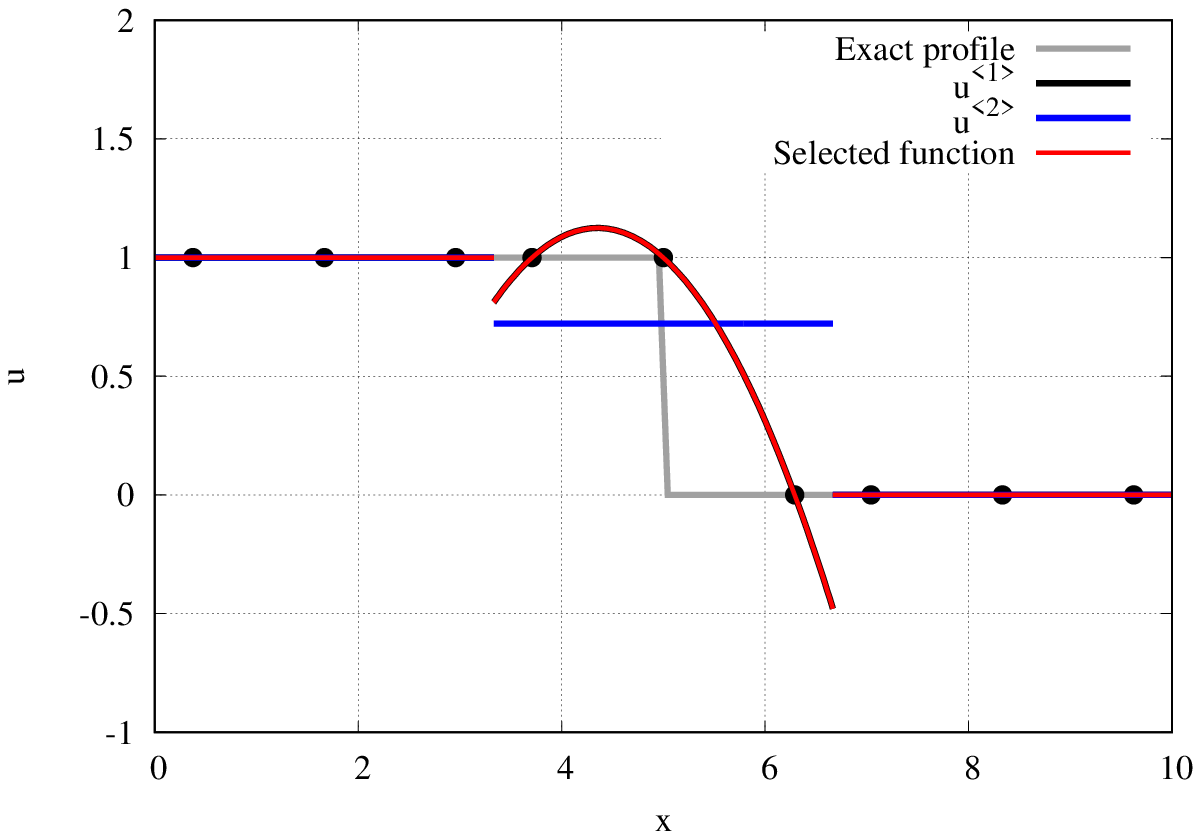}}
\caption{Polynomial selection for the Case 2 solution (discontinuous flow given by Eq.\eqref{eq_Case2}); polynomial degree $K=2$.
The grey lines show the analytical solution;
black lines show the high-order polynomial approximation $u^{<1>}_i$;
blue lines show the stable function $u^{<2>}_i$;
red lines show the selected function based on TVB or BVD criteria.
The computational domain is divided into three cells, where the functions in the first and third cell are fixed to $u^{<1>}_i$
 }\label{fig_Case2_K2}
\end{figure}

\clearpage
\subsection{Linear advection equation}
In this subsection, the BVD algorithm and TVB limiter are examined for a linear advection equation:
\begin{align}
\frac{\del u}{\del t} + \frac{\del u}{\del x}=0.\label{eq_Advec}
\end{align}
The computational domain is set to be $0\leq x\leq 10$;
the number of total degrees of freedom (DoF) is $240$, i.e., the number of total cells is $80$ with $K=2$.
The initial solutions $u(0,x)$ are given as:
\begin{itemize}
\item Case 5: discontinuous flow
\begin{align}
u(0,x)=
\begin{cases}
1.0&\quad (x\leq 1.0),\label{eq_Case5}\\
0.0&\quad ( 1.0<x),
\end{cases}
\end{align}
\item Case 6: convex, smooth, and discontinuous flow 
\begin{align}
u(0,x)&=
\begin{cases}
A\sin(2\pi\omega x)&\quad (x\leq 1.5),\\
1.0+A\sin(2\pi\omega x)&\quad (1.5 < x\leq 3.0),\\
A\sin(2\pi\omega x)&\quad (3.0 < x),
\end{cases}\label{eq_Case6}\\
A&=0.20,\quad \omega=0.20.
\end{align}
\end{itemize}

Figure \ref{fig_Case5_K2} shows the computational results of the case 5 (Eq.\eqref{eq_Case5}) with $K=2$.
The grey-colored lines show the theoretical solution;
the blue and red lines indicate the numerical solution at $t=0.8$ and $8.0$, respectively.
The important features are summarized as follows:
\begin{itemize}
\item
Both of the BVD algorithm and TVB limiter achieve stable computations using the WENO reconstruction in Figs.\ref{fig_Case5_K2}(a)--(c).
However, the BVD case shows more dissipative solution than the conventional TVB limiter cases.
\item
The simple WENO reconstruction in Figs.\ref{fig_Case5_K2}(d)--(f), on the other hand, shows that the BVD algorithm provides sufficiently accurate solution comparable to the TVB limiter cases.
Furthermore, the small oscillation near the discontinuity is removed in the BVD case with maintaining less dissipative solution than the other TVB limiter cases.
\end{itemize}
The computational results of the case 6 (Eq.\eqref{eq_Case6}) with $K=2$ are shown in Fig.\ref{fig_Case6_K2}, where the trend of the BVD algorithm and TVB limiter is almost the same as that discussed in the case 5.

Overall, the BVD framework with the WENO and simple WENO reconstruction effectively works and achieves stable computations without any ad hoc parameters, which would be more reliable than the existent TVB limiter method.
Furthermore, in the present test cases, the BVD algorithm with the simple WENO reconstruction would provides better result than the conventional TVB limiter cases in the FR framework (possiblly in the nodal-type DG scheme as well).
This is also supported by the higher-order cases ($K=3$) shown in Appendix B.
%
%

\renewcommand{\datahead}{bv_data_cfl0p1_lhs3_rhs1_cor0_sp0_init}
\renewcommand{\timestep}{00001440}
\renewcommand{\wid}{0.3275} 
\begin{figure}[!htbp]
    \centering
    \setcounter{subfigure}{0}
    \subfigure[][{\scriptsize BVD/WENO}]{\includegraphics[width=\wid\textwidth,clip]{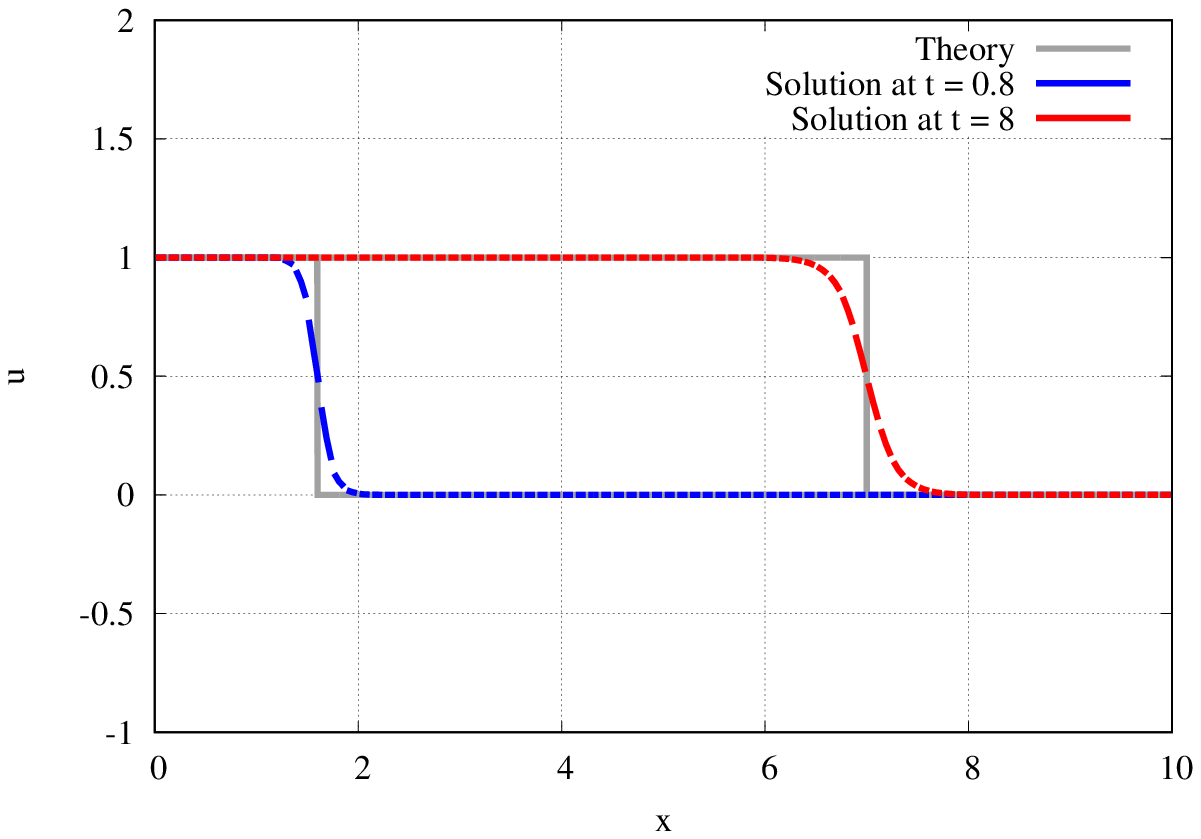}}
    \subfigure[][{\scriptsize TVB/WENO $(M=0.0)$}]{\includegraphics[width=\wid\textwidth,clip]{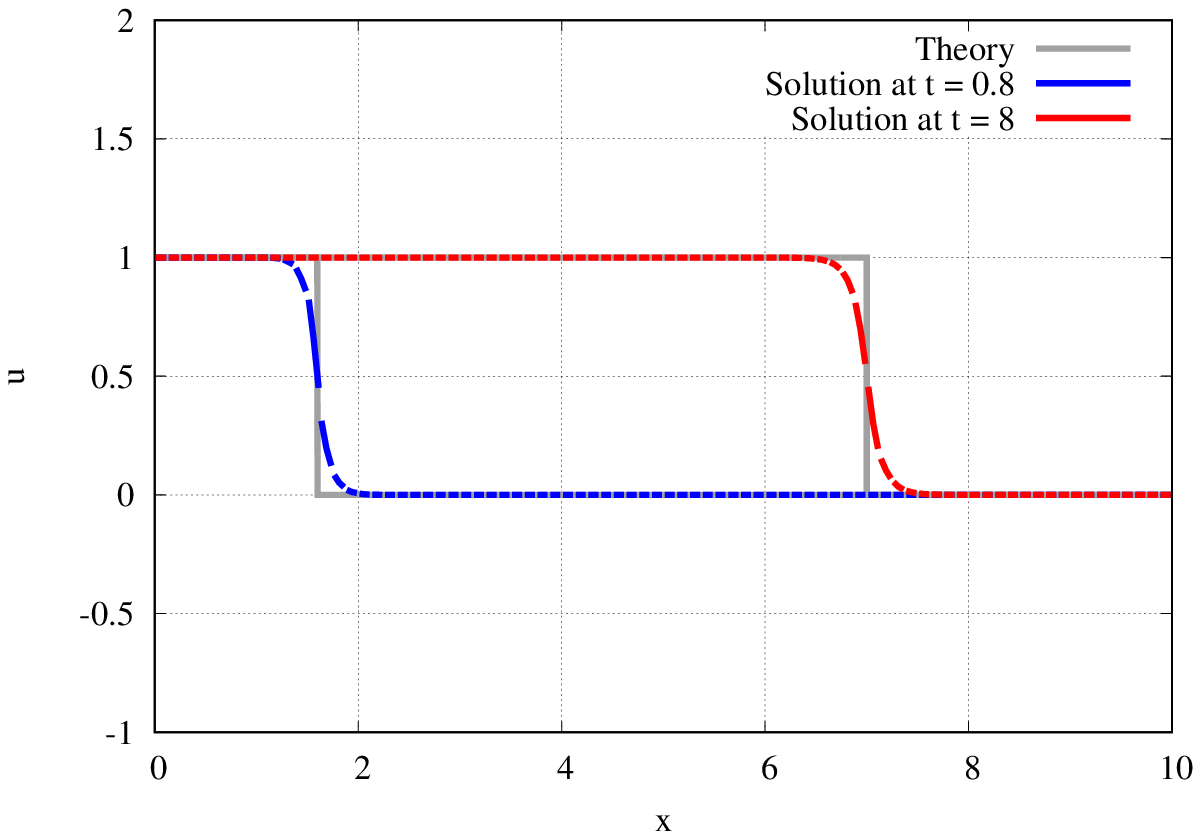}}
    \subfigure[][{\scriptsize TVB/WENO $(M=200)$}]{\includegraphics[width=\wid\textwidth,clip]{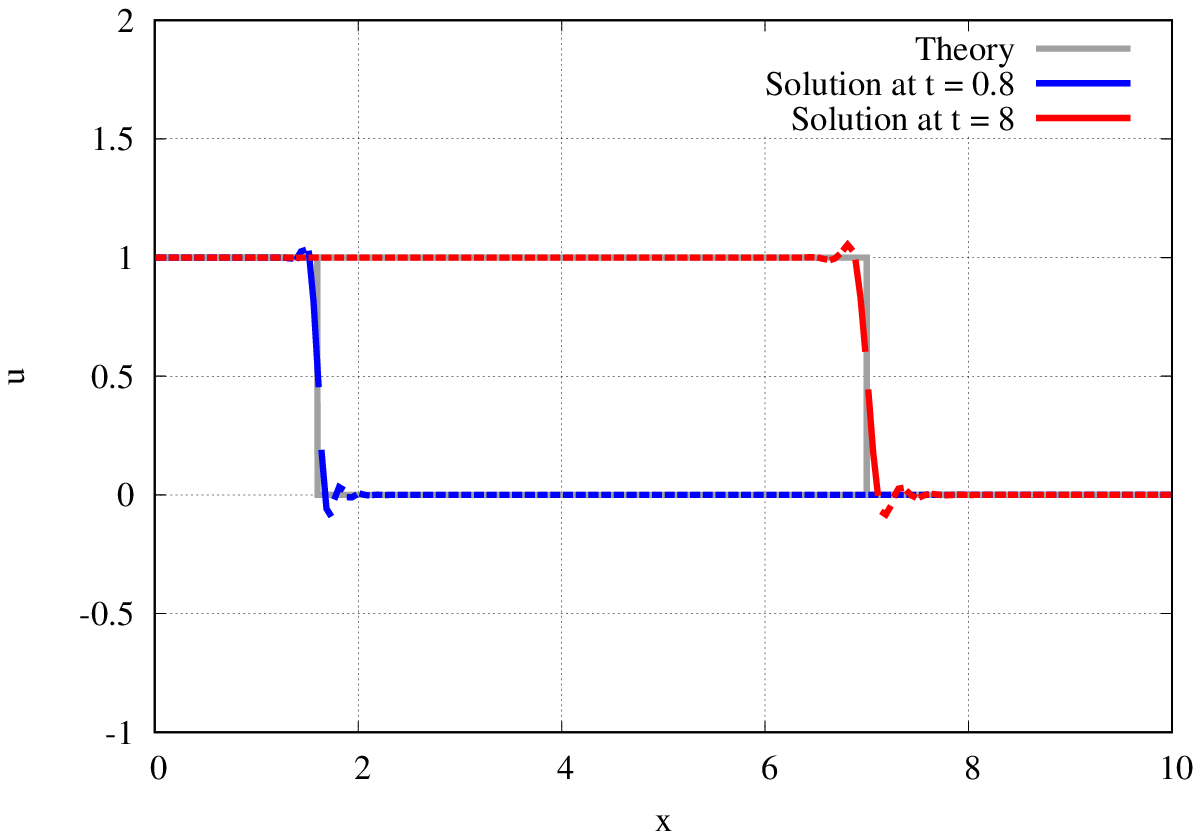}}\\
    \subfigure[][{\scriptsize BVD/SWENO}]{\includegraphics[width=\wid\textwidth,clip]{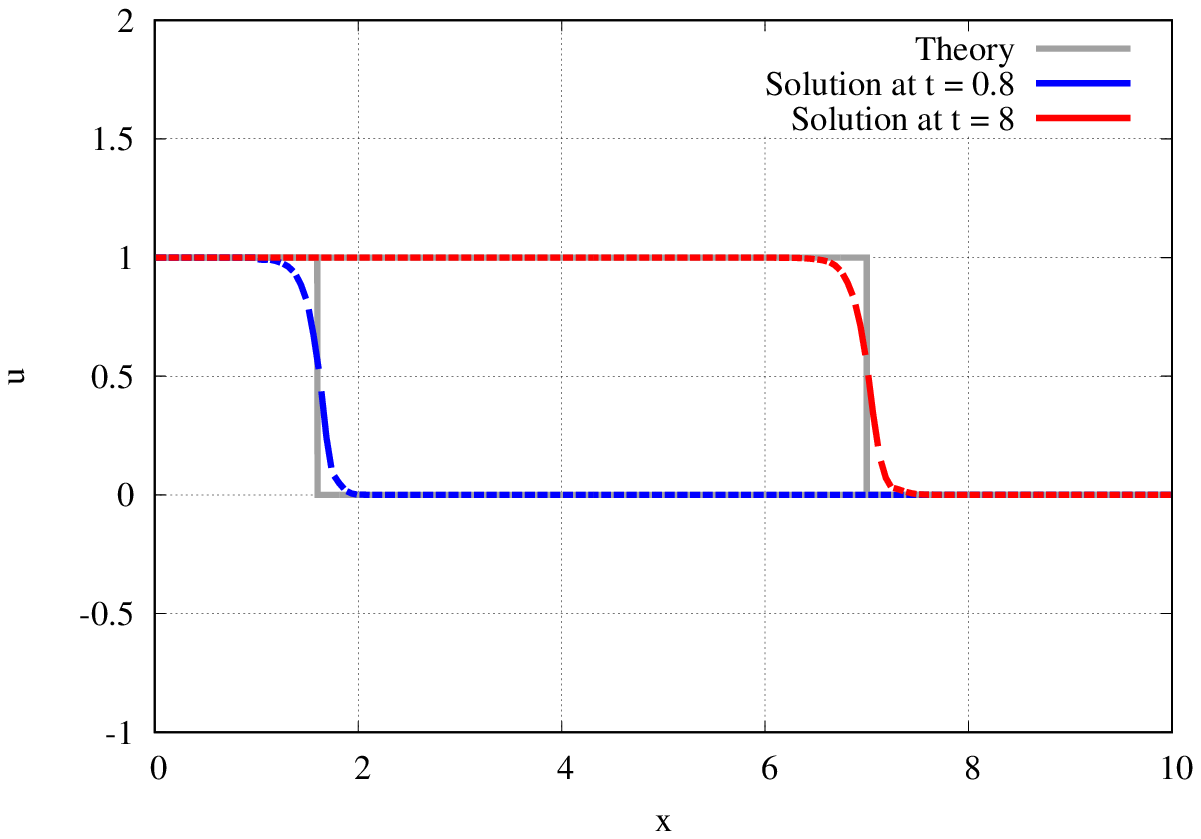}}
    \subfigure[][{\scriptsize TVB/SWENO $(M=0.0)$}]{\includegraphics[width=\wid\textwidth,clip]{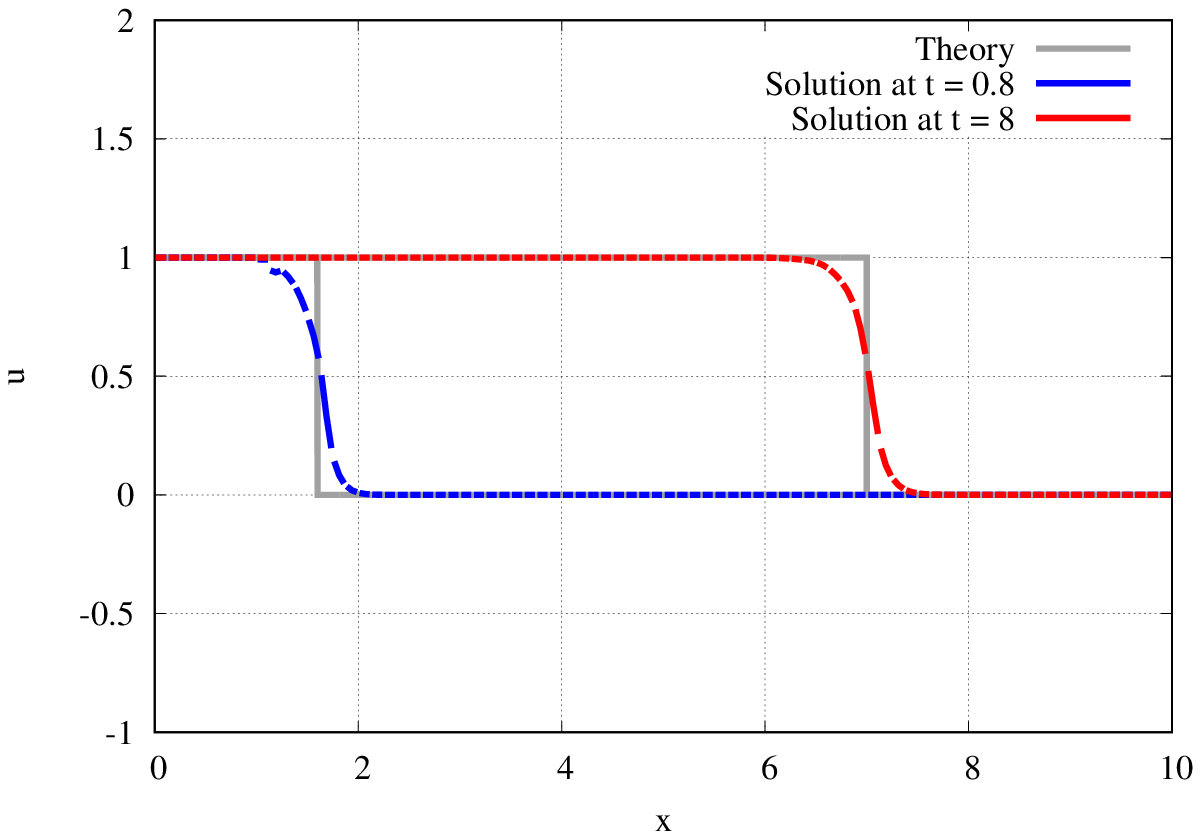}}
    \subfigure[][{\scriptsize TVB/SWENO $(M=200)$}]{\includegraphics[width=\wid\textwidth,clip]{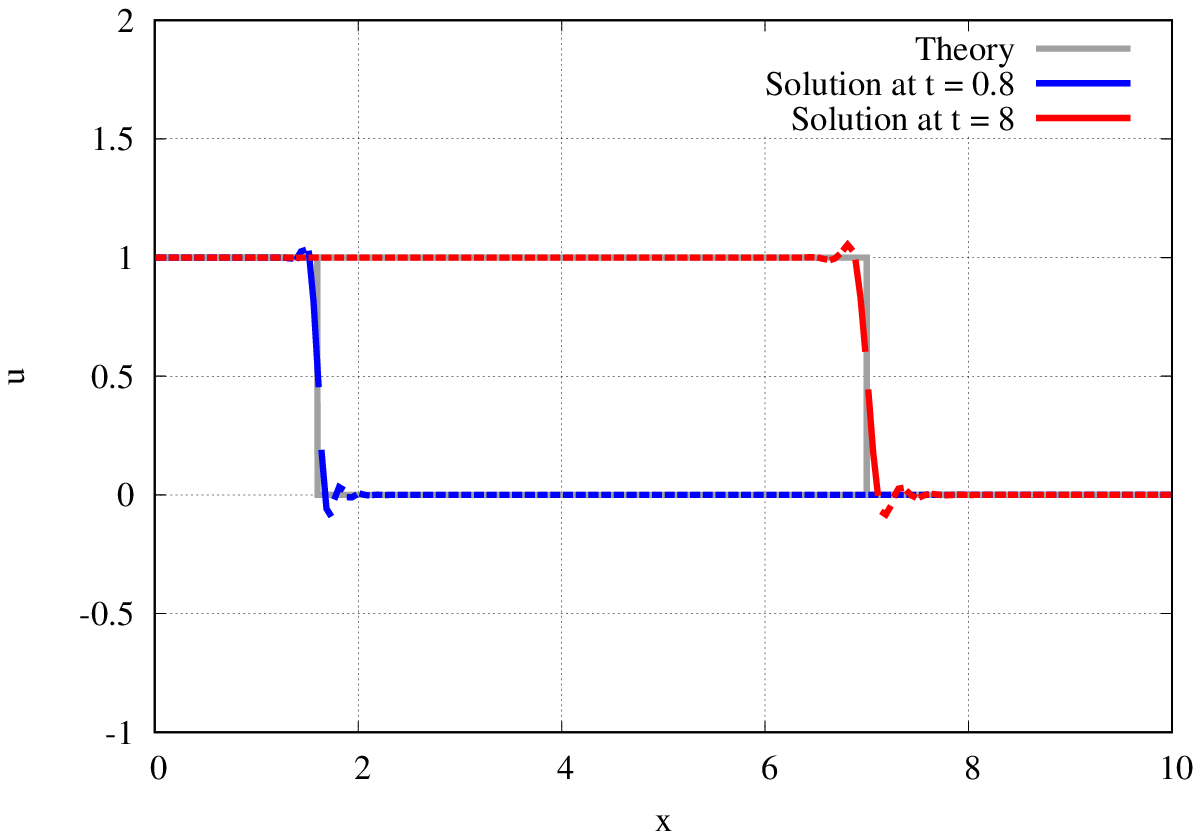}}
\caption{Computational results of a linear advection equation for the case 5 initial condition (discontinuous flow given by Eq.\eqref{eq_Case5}); the total DoF is $240$ with a polynomial degree $K=2$.
The grey lines show the theoritical solution;
blue lines show the solution $u(x,t)$ at $t=0.8$;
red lines show the solution $u(x,t)$ at $t=8.0$.
 }\label{fig_Case5_K2}
\end{figure}
\begin{figure}[!htbp]
    \centering
    \setcounter{subfigure}{0}
    \subfigure[][{\scriptsize BVD/WENO}]{\includegraphics[width=\wid\textwidth,clip]{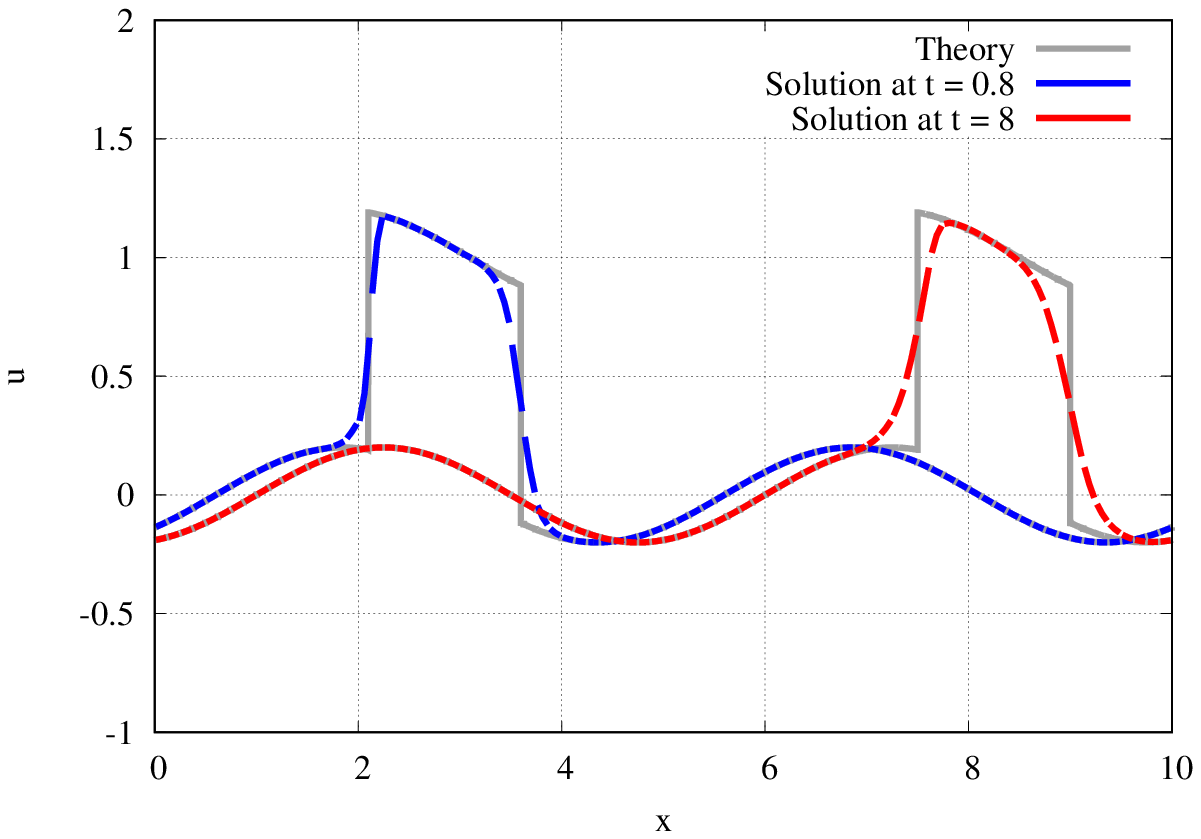}}
    \subfigure[][{\scriptsize TVB/WENO $(M=0.0)$}]{\includegraphics[width=\wid\textwidth,clip]{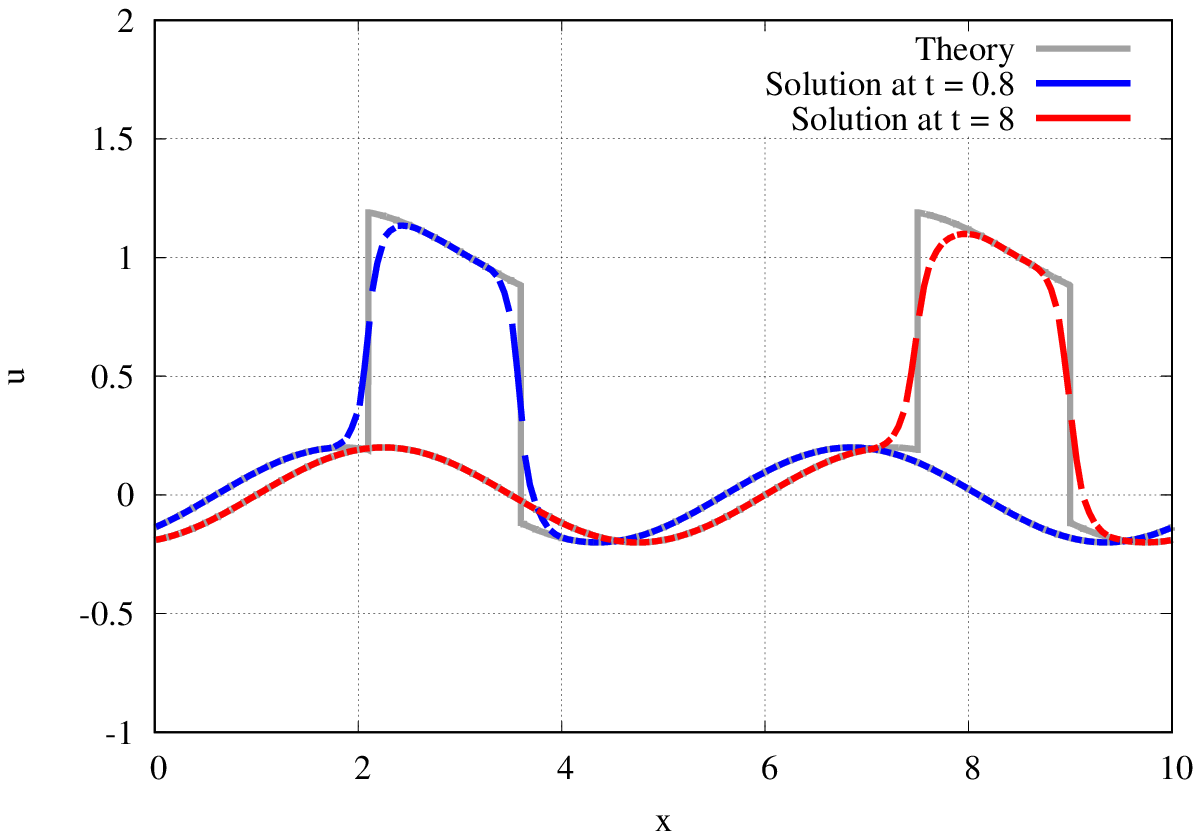}}
    \subfigure[][{\scriptsize TVB/WENO $(M=200)$}]{\includegraphics[width=\wid\textwidth,clip]{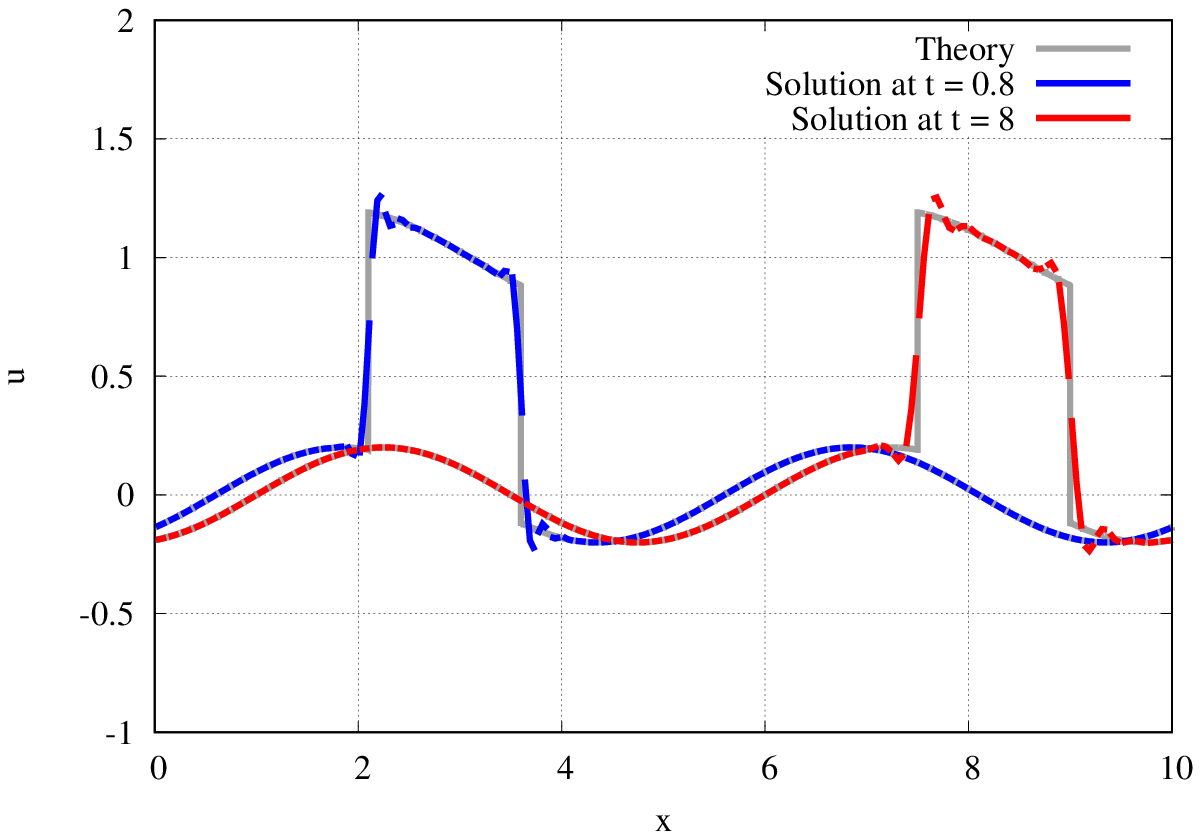}}\\
    \subfigure[][{\scriptsize BVD/SWENO}]{\includegraphics[width=\wid\textwidth,clip]{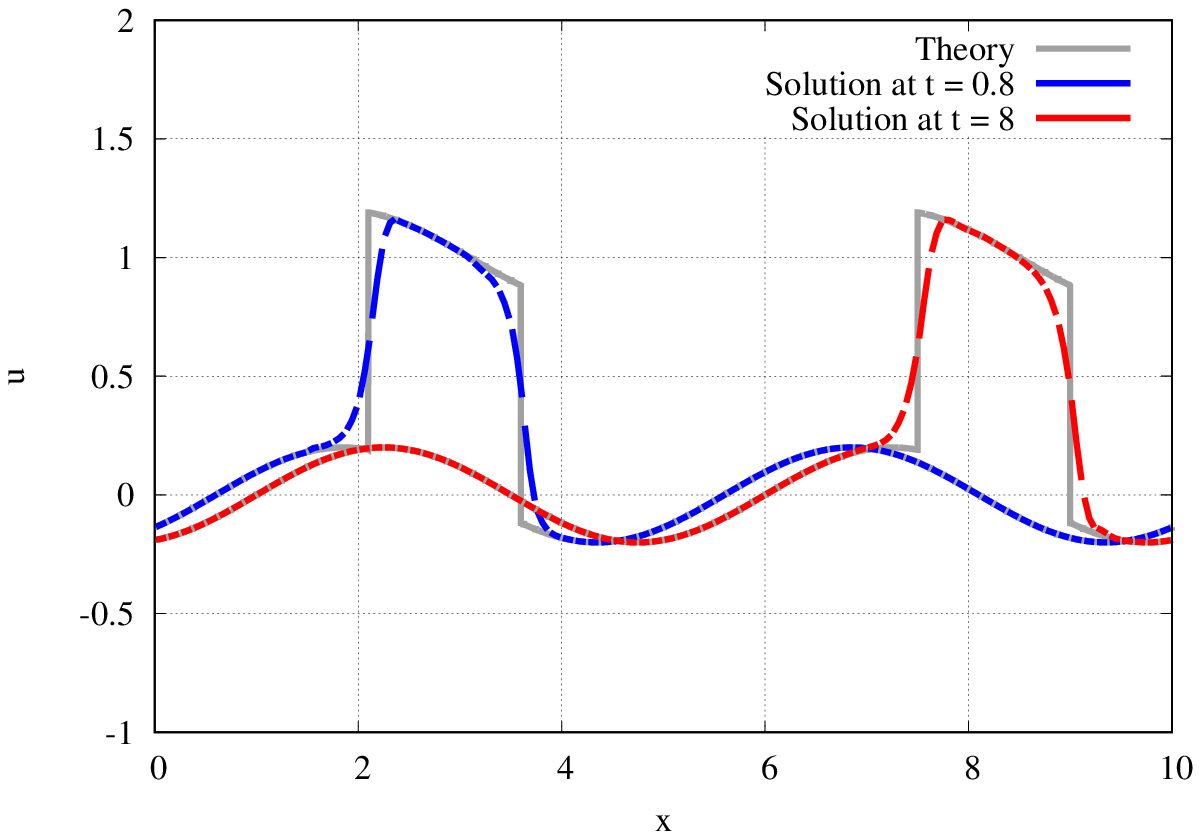}}
    \subfigure[][{\scriptsize TVB/SWENO $(M=0.0)$}]{\includegraphics[width=\wid\textwidth,clip]{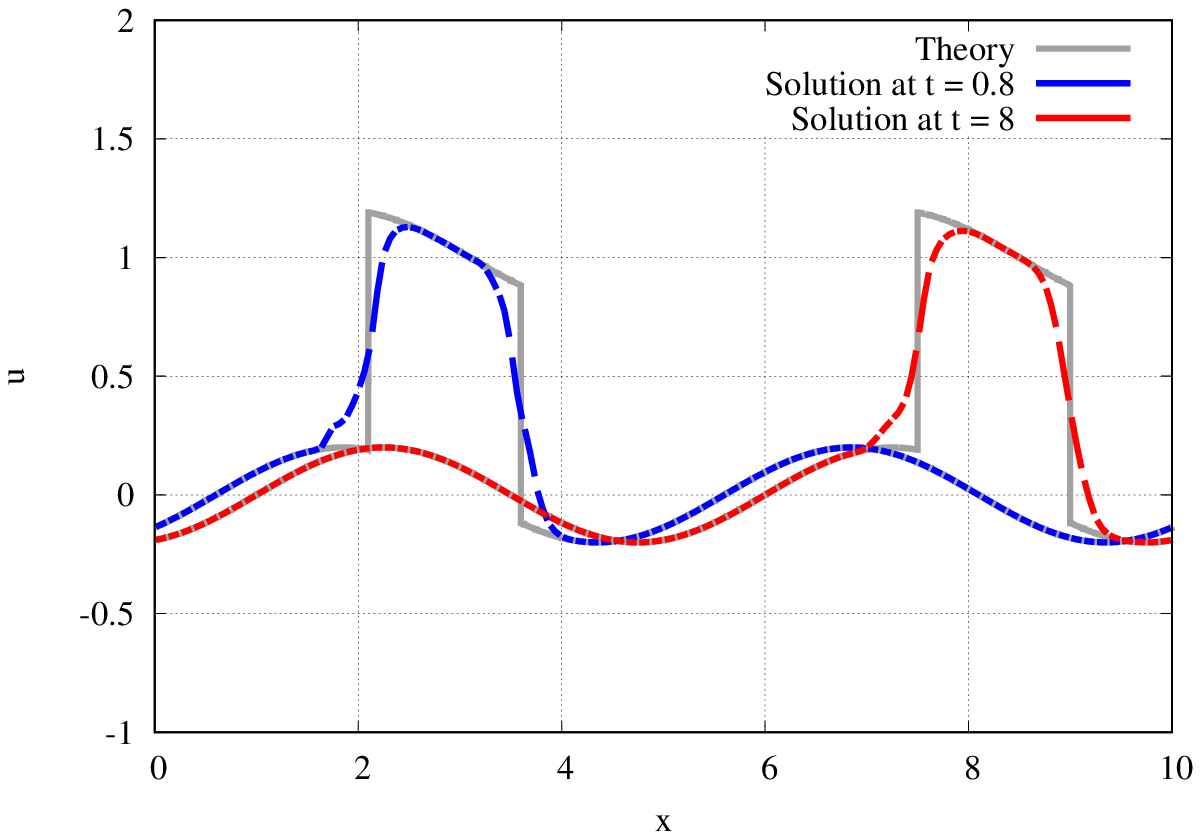}}
    \subfigure[][{\scriptsize TVB/SWENO $(M=200)$}]{\includegraphics[width=\wid\textwidth,clip]{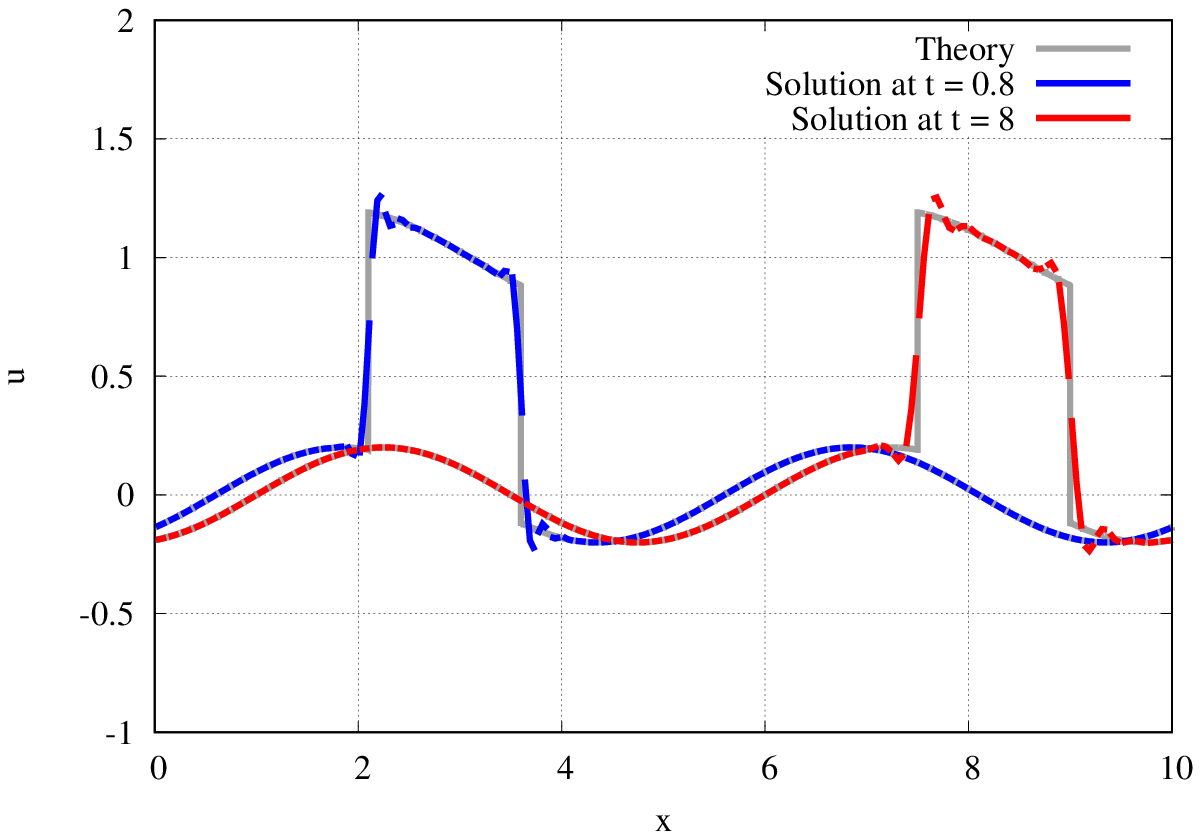}}
\caption{Computational results of a linear advection equation for the case 6 initial condition (convex, smooth, and discontinuous flow given by Eq.\eqref{eq_Case6}); the total DoF is $240$ with a polynomial degree $K=2$.
The grey lines show the theoritical solution;
blue lines show the solution $u(x,t)$ at $t=0.8$;
red lines show the solution $u(x,t)$ at $t=8.0$.
 }\label{fig_Case6_K2}
\end{figure}

\section{Summary}
This paper represents an implementation of a new algorithm for appropriate reconstructions, so-called the BVD algorithm, to the hyperbolic conservation laws in the FR framework. 
The applicability of the BVD algorithm to the WENO methodology was numerically examined for a linear advection equation as well as the selection of appropriate reconstructions for a given solution profile.
\begin{itemize}
\item
In the WENO reconstruction methodology, the polynomial selection by the BVD algorithm is not always the same as those by the conventional TVB limiter depending on the ad hoc TVB parameter.
\item
Both of the BVD algorithm and TVB limiter attain stable computations of a linear advection equation with a discontinuous initial conditions.
In the present test cases, the BVD algorithm with the simple WENO reconstruction provides more stable and accurate solution compared to the conventional TVB limiter cases;
however, the reconstruction based on the WENO method results in more dissipative solutions using the BVD algorithm than TVB limiter one.
\item
Overall, the BVD algorithm realizes a completely paramter-free selection unlike the conventional indicator based on the TVB limiter, which would possiblly provide more reliable results and offer a radically new criterion.
\end{itemize}
The present trend would be expected also in the system equations such as the Euler equations although its implementation should be carefully discussed.
Finally, the present work is limited to third or lower order polynomials, leaving the implementations for higher-order schemes a future work.
Furthermore, for the simplicity, this note only focuses on a linear advection equation, and the nonlinear and system equations will be discussed in the next paper.

\bibliographystyle{elsart-num-sort}
\bibliography{flab}

\clearpage
\setcounter{section}{0}
\renewcommand{\thesection}{\Alph{section}}
\section{Polynomial selection by the TVB and BVD criteria: advanced cases}
In this section, the advanced cases in the numerical test for polynomial selection are presented as follows:
\begin{itemize}
\item Case 3: smooth and discontinuous flow
\begin{align}
u&=
\begin{cases}
1.0+A\sin(2\pi\omega x)&\quad (x\leq 5.0),\\
A\sin(2\pi\omega x)&\quad ( 5.0<x),
\end{cases}\label{eq_Case3}\\
A&=0.20,\quad \omega=0.25,
\end{align}
\item Case 4: convex, smooth, and discontinuous flow 
\begin{align}
u&=
\begin{cases}
A\sin(2\pi\omega x)&\quad (x\leq 3.5),\\
1.0+A\sin(2\pi\omega x)&\quad (3.5 < x\leq 5.0),\\
A\sin(2\pi\omega x)&\quad (5.0 < x),
\end{cases}\label{eq_Case4}\\
A&=0.20,\quad \omega=0.25.
\end{align}
\end{itemize}
These cases comprise the combination of the smooth and discontinuous solutions.
The higher-order approximation $K=3$ is additionally presented in Fig.\ref{fig_Case3_K3} for the simple WENO reconstruction in the case 3 (smooth and discontinuous flow: Eq.\eqref{eq_Case3}).
\renewcommand{\datahead}{bv_data_cfl0p1_lhs3_rhs1_cor1_sp0_init}
\renewcommand{\timestep}{00000001}
\renewcommand{\wid}{0.3275} 
\begin{figure}[!htbp]
    \centering
    \setcounter{subfigure}{0}
    \subfigure[][{\scriptsize BVD/WENO}]{\includegraphics[width=\wid\textwidth,clip]{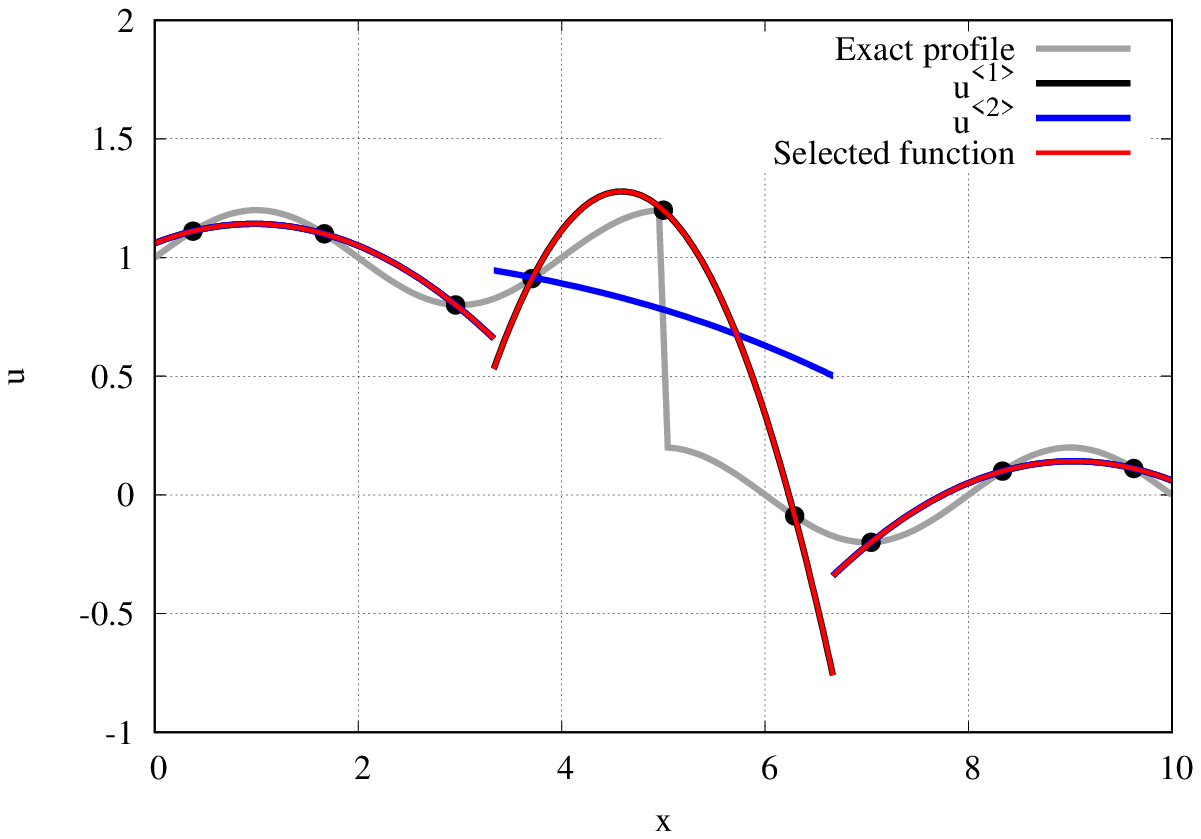}}
    \subfigure[][{\scriptsize TVB/WENO $(M=0.0)$}]{\includegraphics[width=\wid\textwidth,clip]{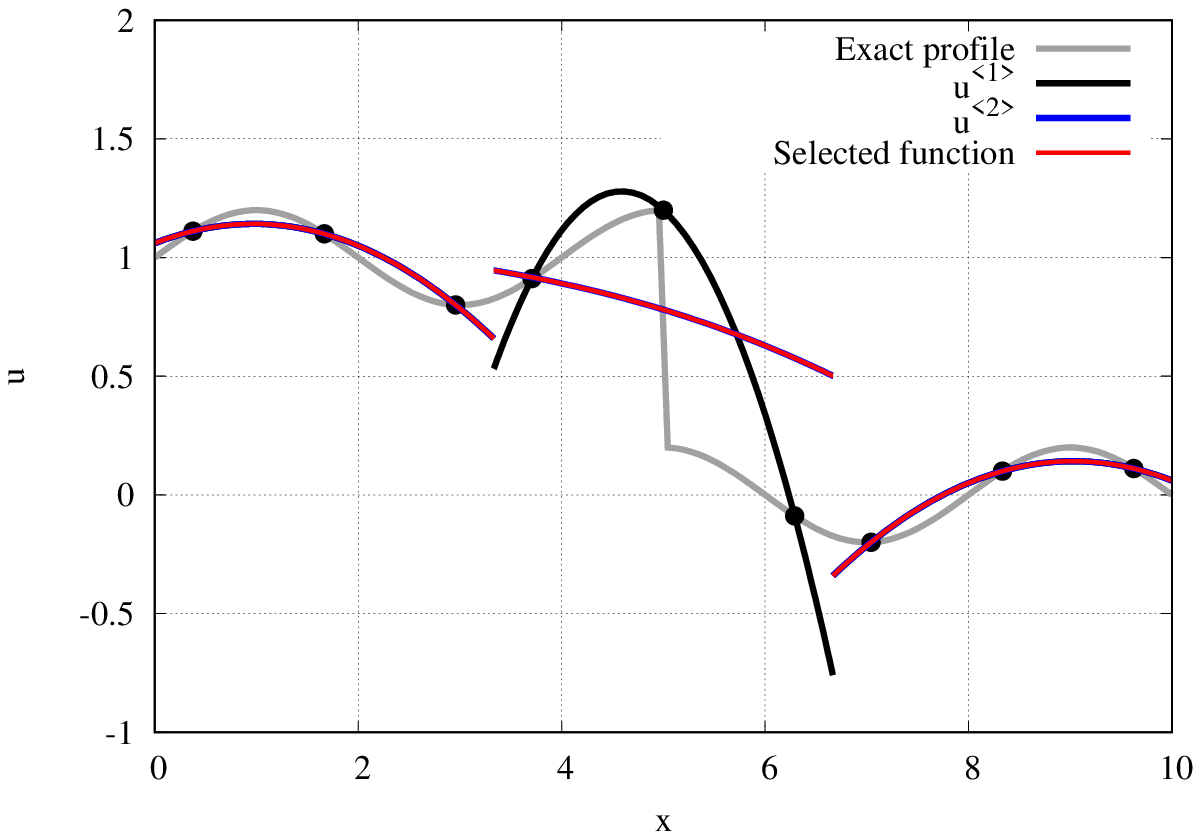}}
    \subfigure[][{\scriptsize TVB/WENO $(M=200)$}]{\includegraphics[width=\wid\textwidth,clip]{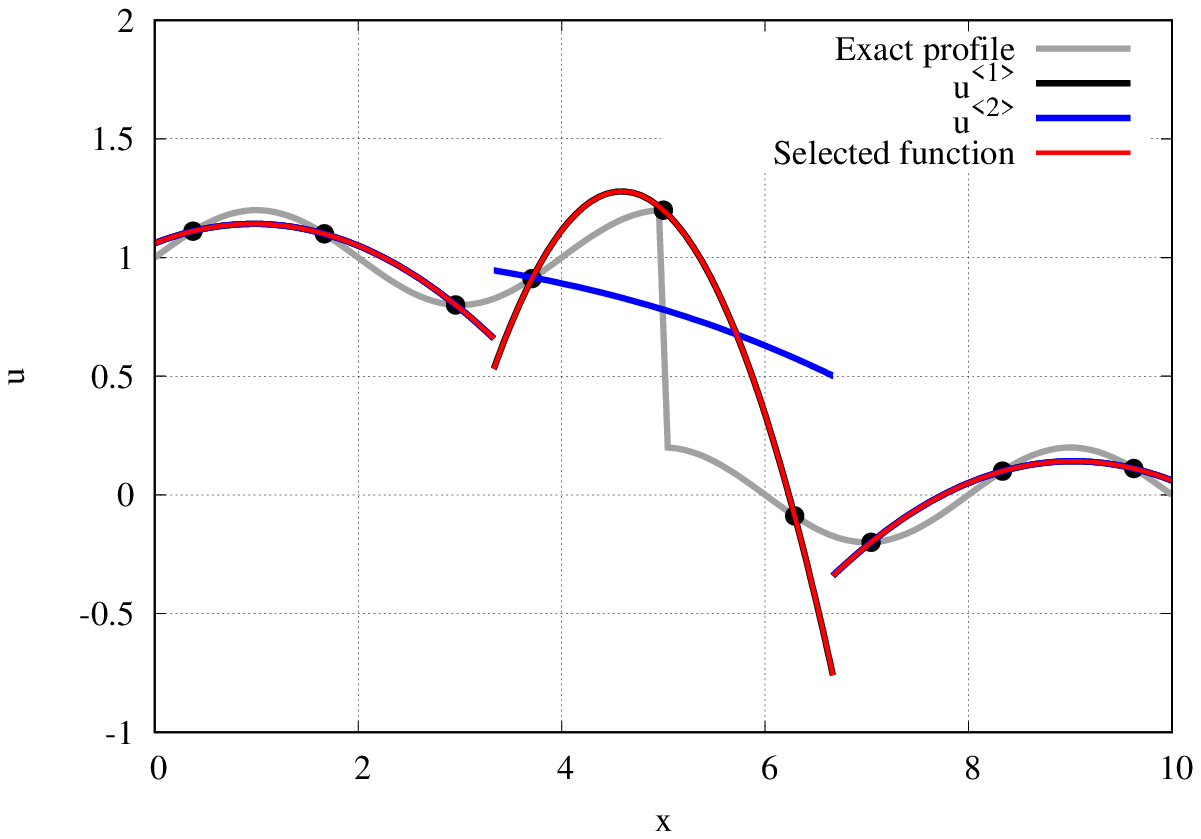}}\\
    \subfigure[][{\scriptsize BVD/SWENO}]{\includegraphics[width=\wid\textwidth,clip]{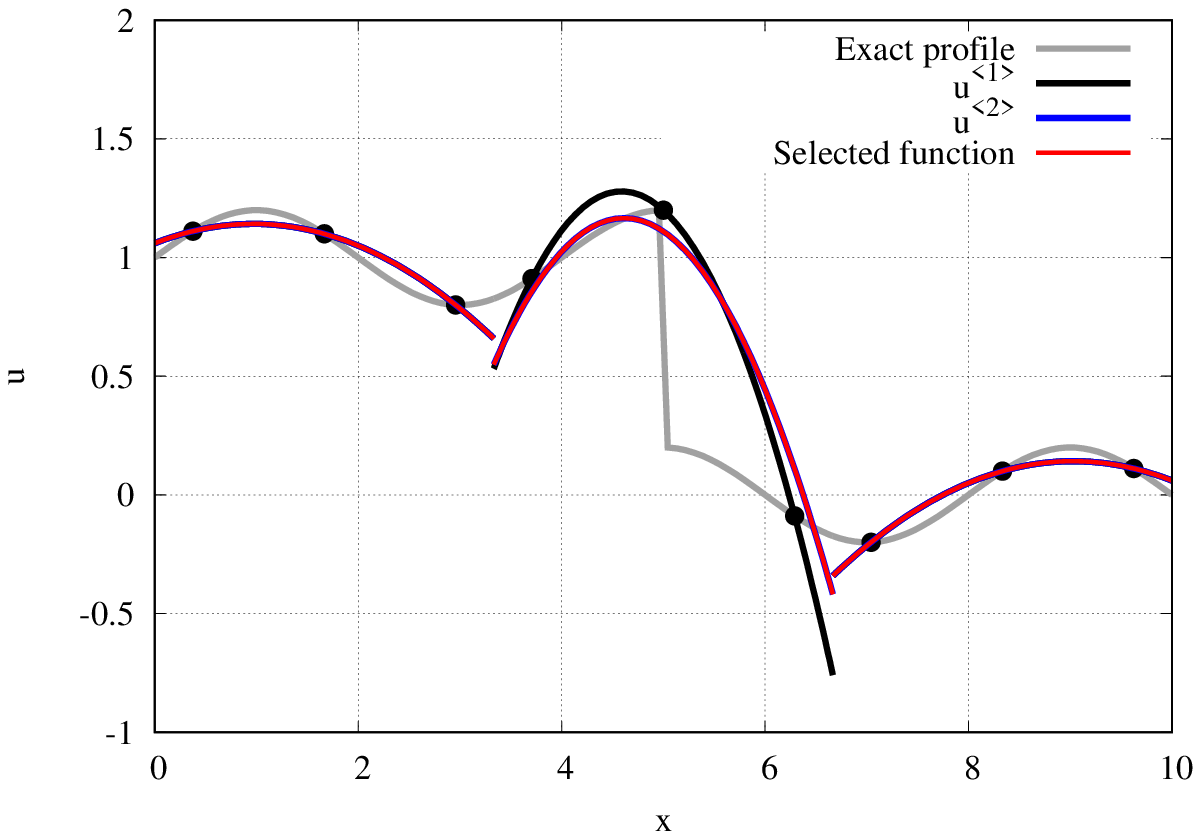}}
    \subfigure[][{\scriptsize TVB/SWENO $(M=0.0)$}]{\includegraphics[width=\wid\textwidth,clip]{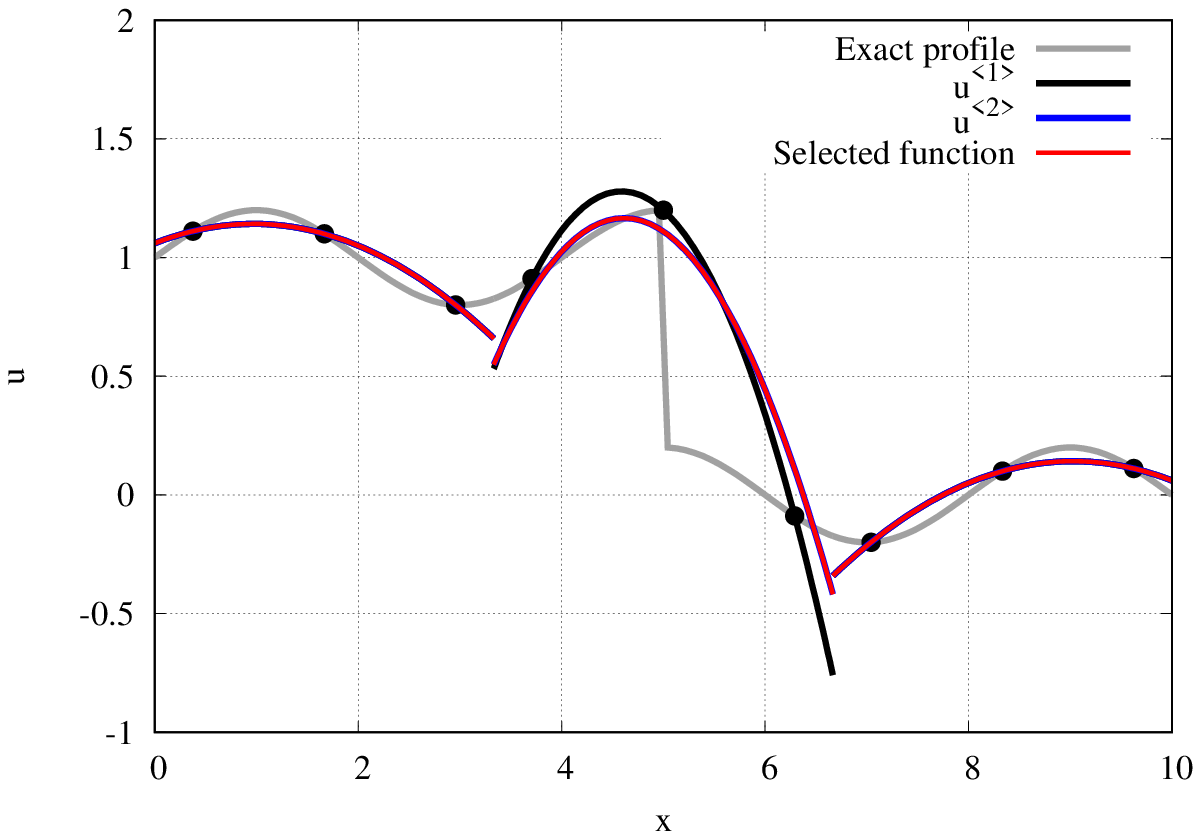}}
    \subfigure[][{\scriptsize TVB/SWENO $(M=200)$}]{\includegraphics[width=\wid\textwidth,clip]{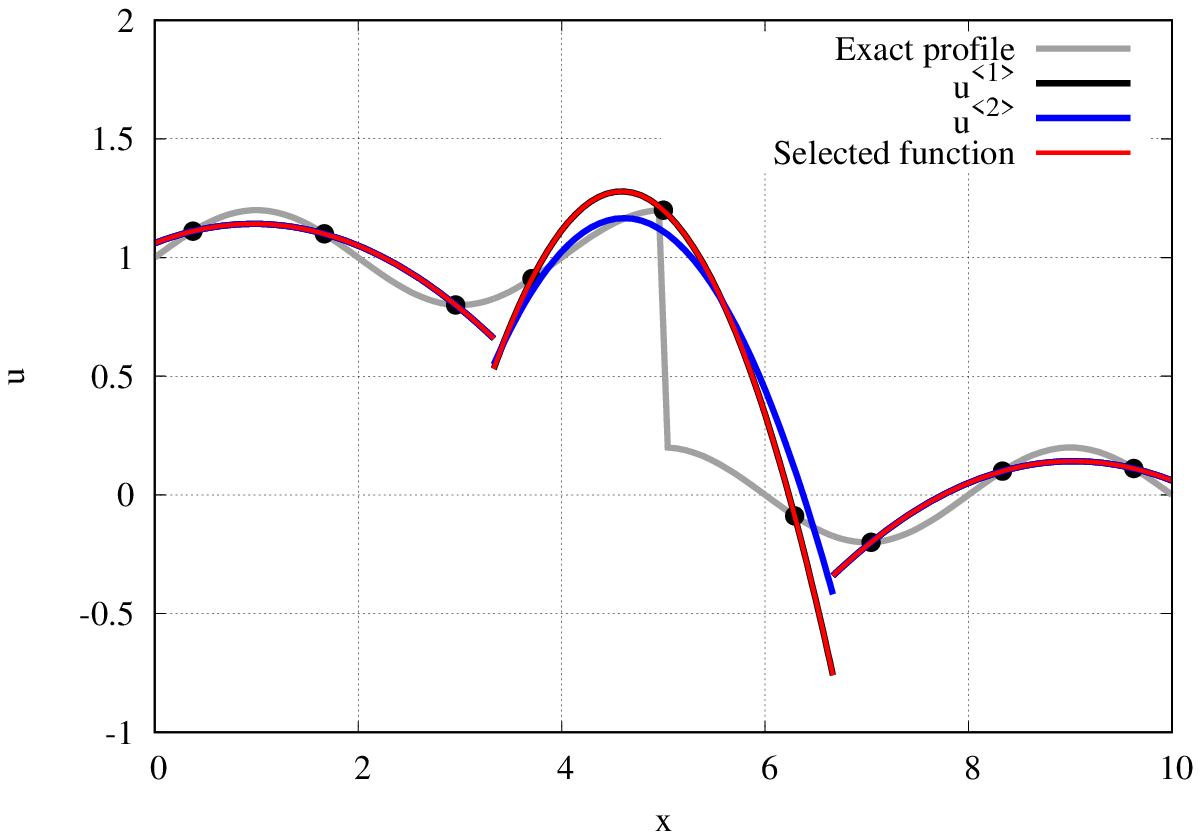}}
\caption{Polynomial selection for the Case 3 solution (smooth and discontinuous flow given by Eq.\eqref{eq_Case3}); polynomial degree $K=2$.
The grey lines show the analytical solution;
black lines show the high-order polynomial approximation $u^{<1>}_i$;
blue lines show the stable function $u^{<2>}_i$;
red lines show the selected function based on TVB or BVD criteria.
The computational domain is divided into three cells, where the functions in the first and third cell are fixed to $u^{<1>}_i$
 }\label{fig_Case3_K2}
\end{figure}
\begin{figure}[!htbp]
    \centering
    \setcounter{subfigure}{0}
    \subfigure[][{\scriptsize BVD/WENO}]{\includegraphics[width=\wid\textwidth,clip]{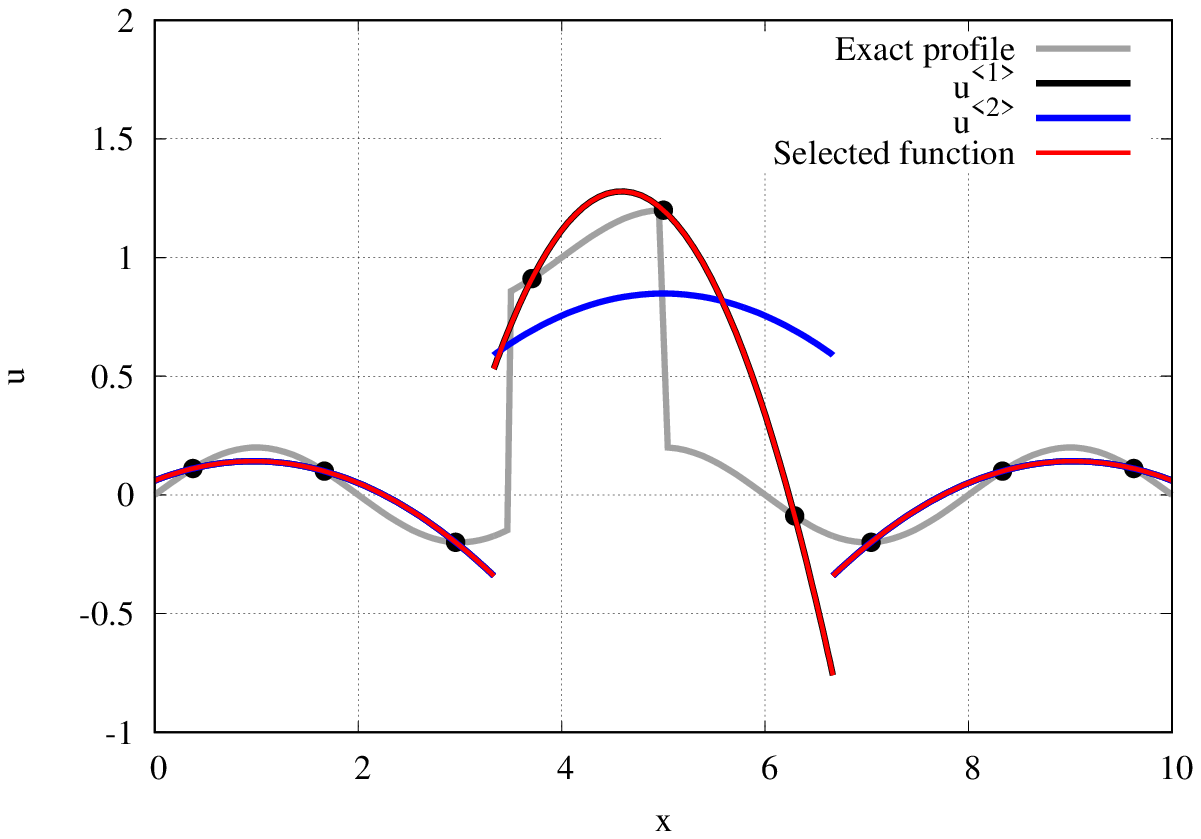}}
    \subfigure[][{\scriptsize TVB/WENO $(M=0.0)$}]{\includegraphics[width=\wid\textwidth,clip]{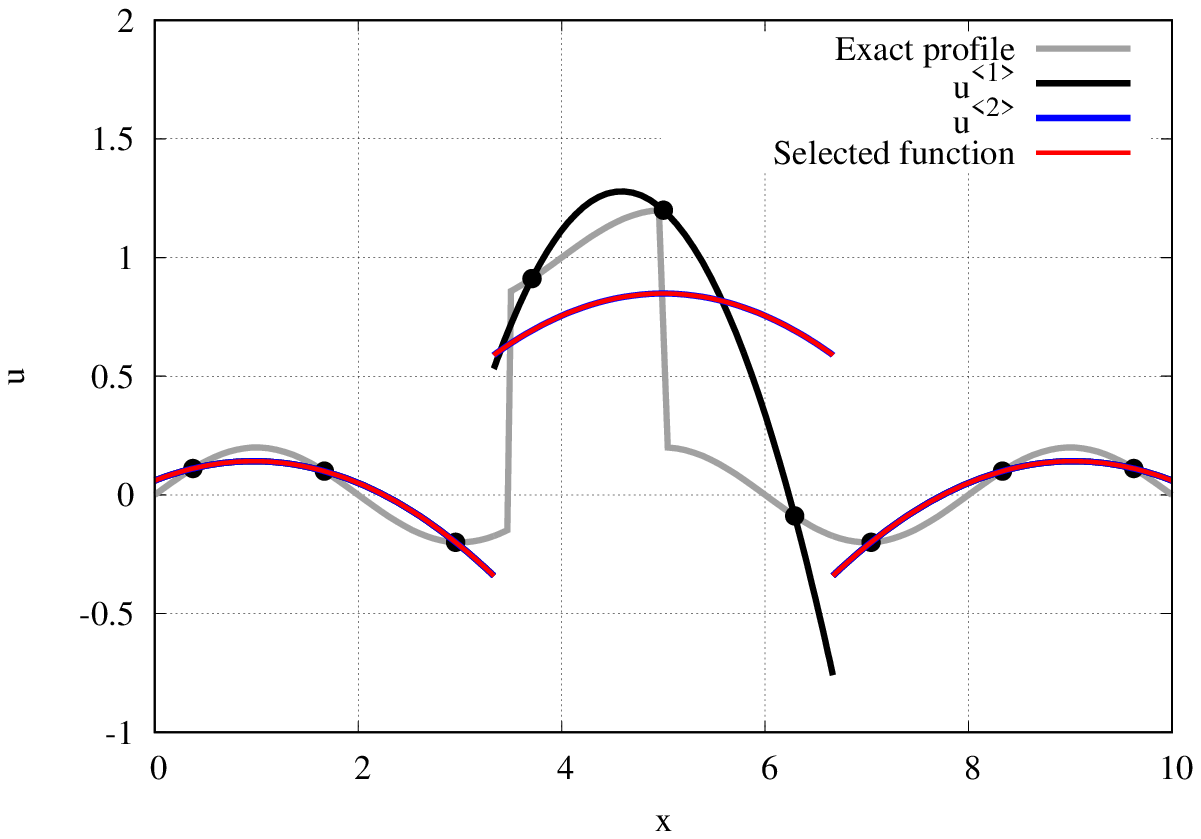}}
    \subfigure[][{\scriptsize TVB/WENO $(M=200)$}]{\includegraphics[width=\wid\textwidth,clip]{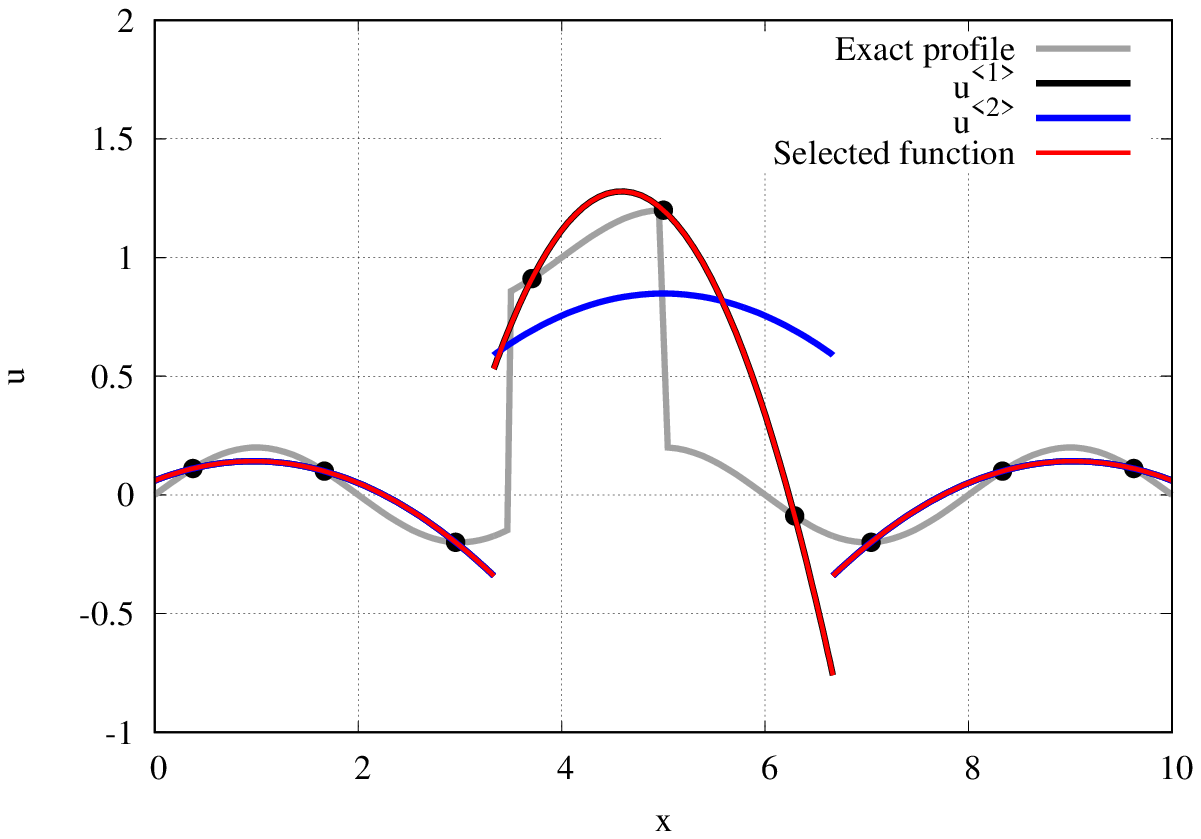}}\\
    \subfigure[][{\scriptsize BVD/SWENO}]{\includegraphics[width=\wid\textwidth,clip]{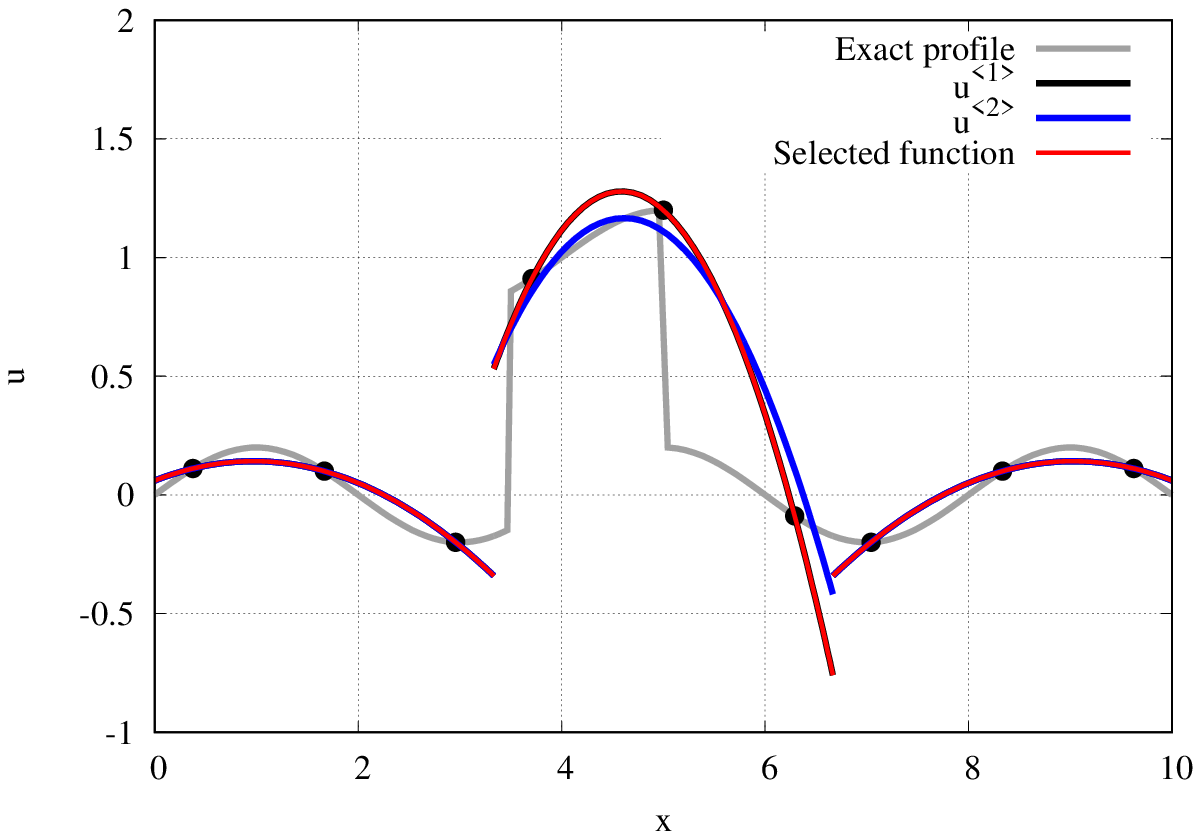}}
    \subfigure[][{\scriptsize TVB/SWENO $(M=0.0)$}]{\includegraphics[width=\wid\textwidth,clip]{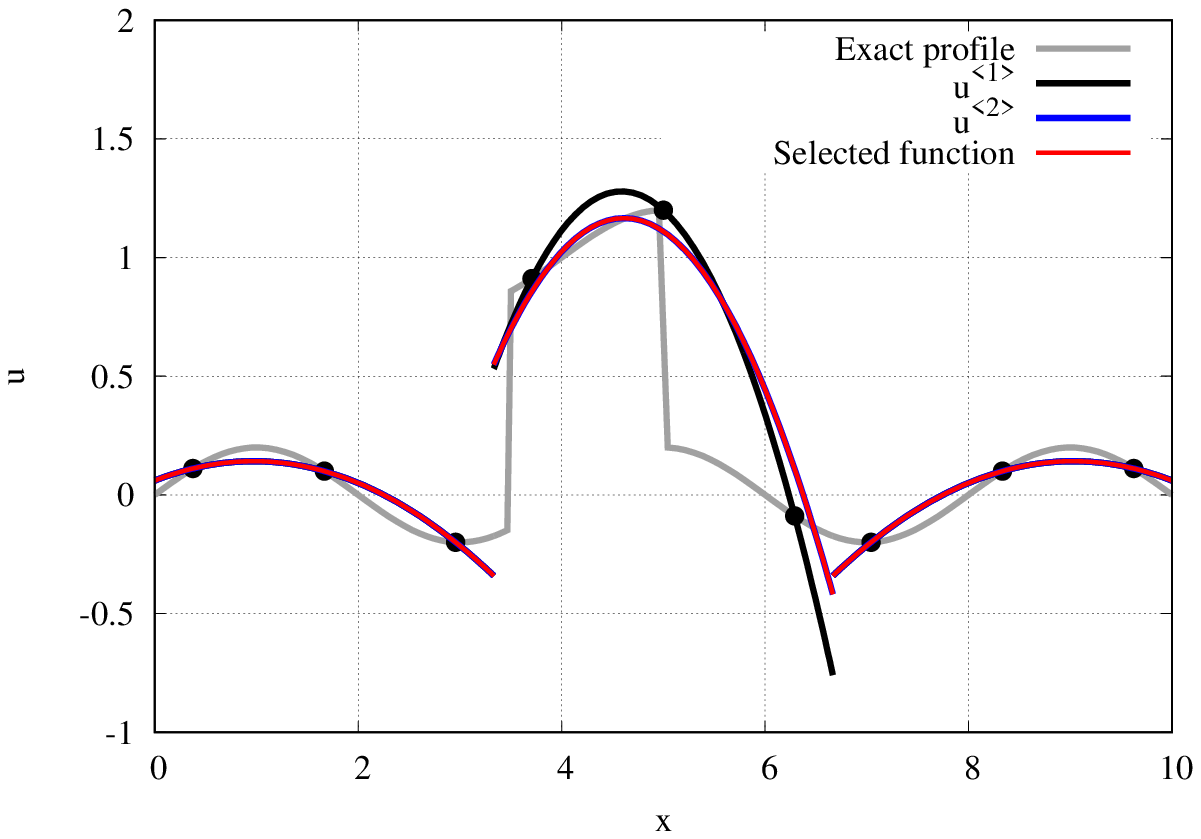}}
    \subfigure[][{\scriptsize TVB/SWENO $(M=200)$}]{\includegraphics[width=\wid\textwidth,clip]{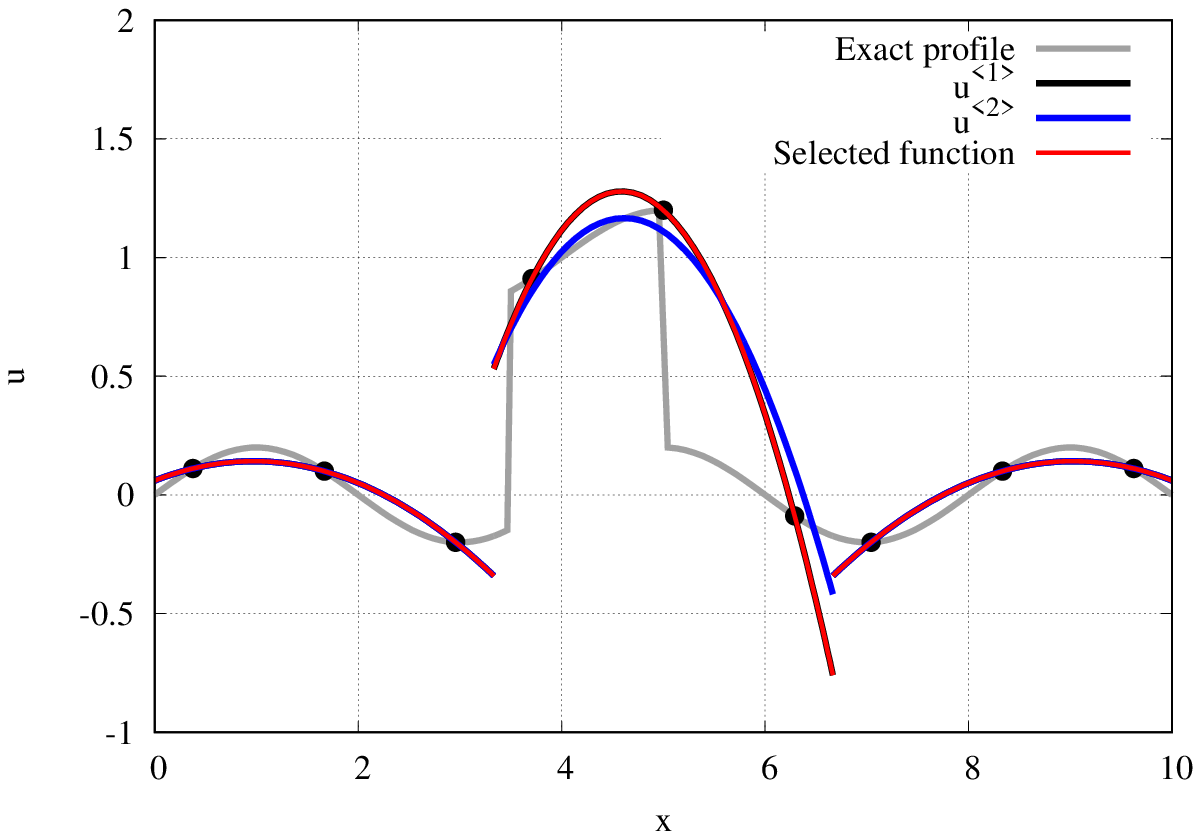}}
\caption{Polynomial selection for the Case 4 solution (convex, smooth, and discontinuous flow given by Eq.\eqref{eq_Case3}); polynomial degree $K=2$.
The grey lines show the analytical solution;
black lines show the high-order polynomial approximation $u^{<1>}_i$;
blue lines show the stable function $u^{<2>}_i$;
red lines show the selected function based on TVB or BVD criteria.
The computational domain is divided into three cells, where the functions in the first and third cell are fixed to $u^{<1>}_i$
 }\label{fig_Case4_K2}
\end{figure}
\begin{figure}[!htbp]
    \centering
    \setcounter{subfigure}{0}
    \subfigure[][{\scriptsize BVD/SWENO}]{\includegraphics[width=\wid\textwidth,clip]{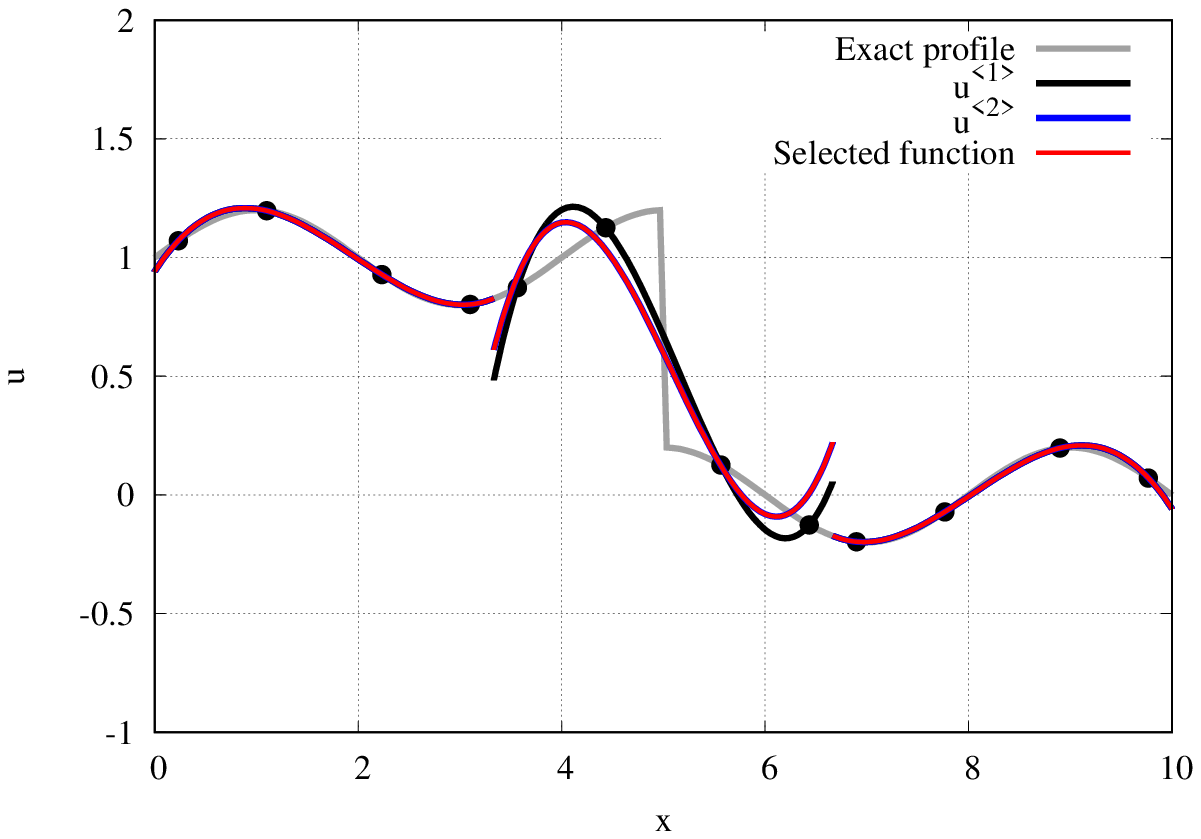}}
    \subfigure[][{\scriptsize TVB/SWENO $(M=0.0)$}]{\includegraphics[width=\wid\textwidth,clip]{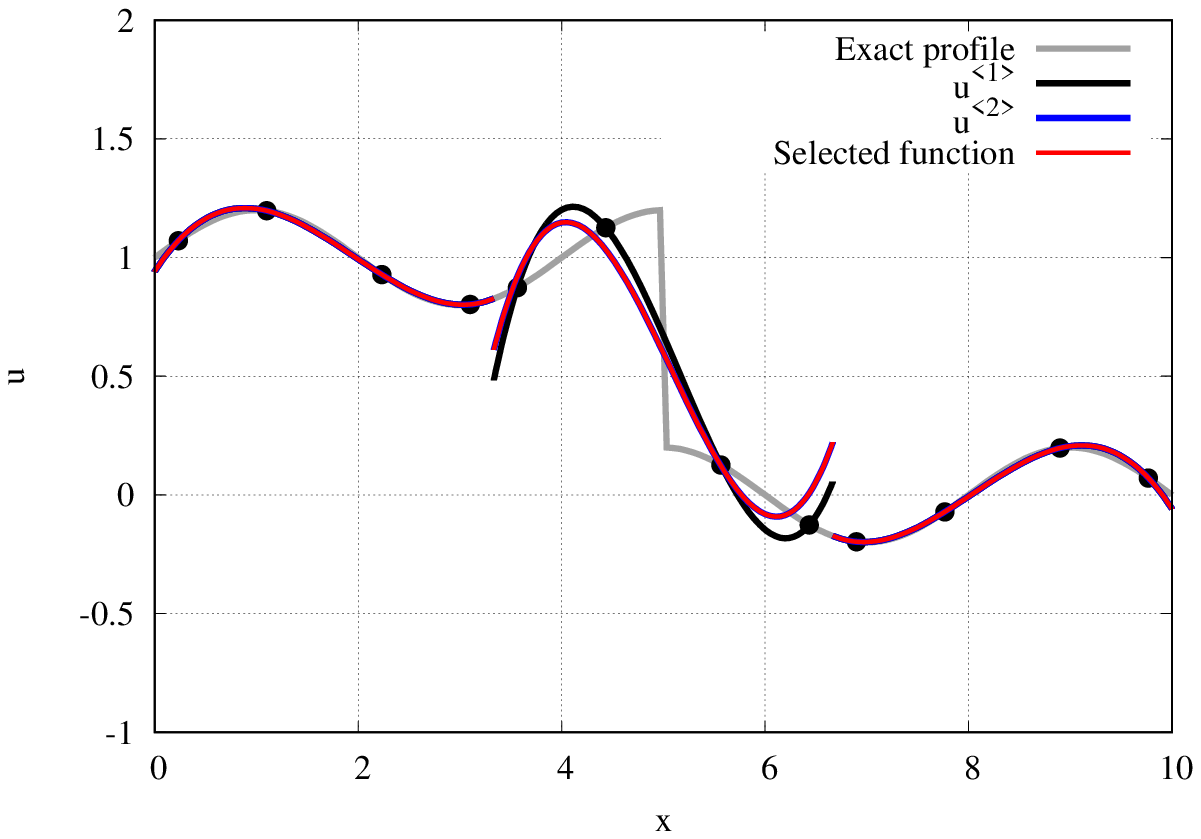}}
    \subfigure[][{\scriptsize TVB/SWENO $(M=200)$}]{\includegraphics[width=\wid\textwidth,clip]{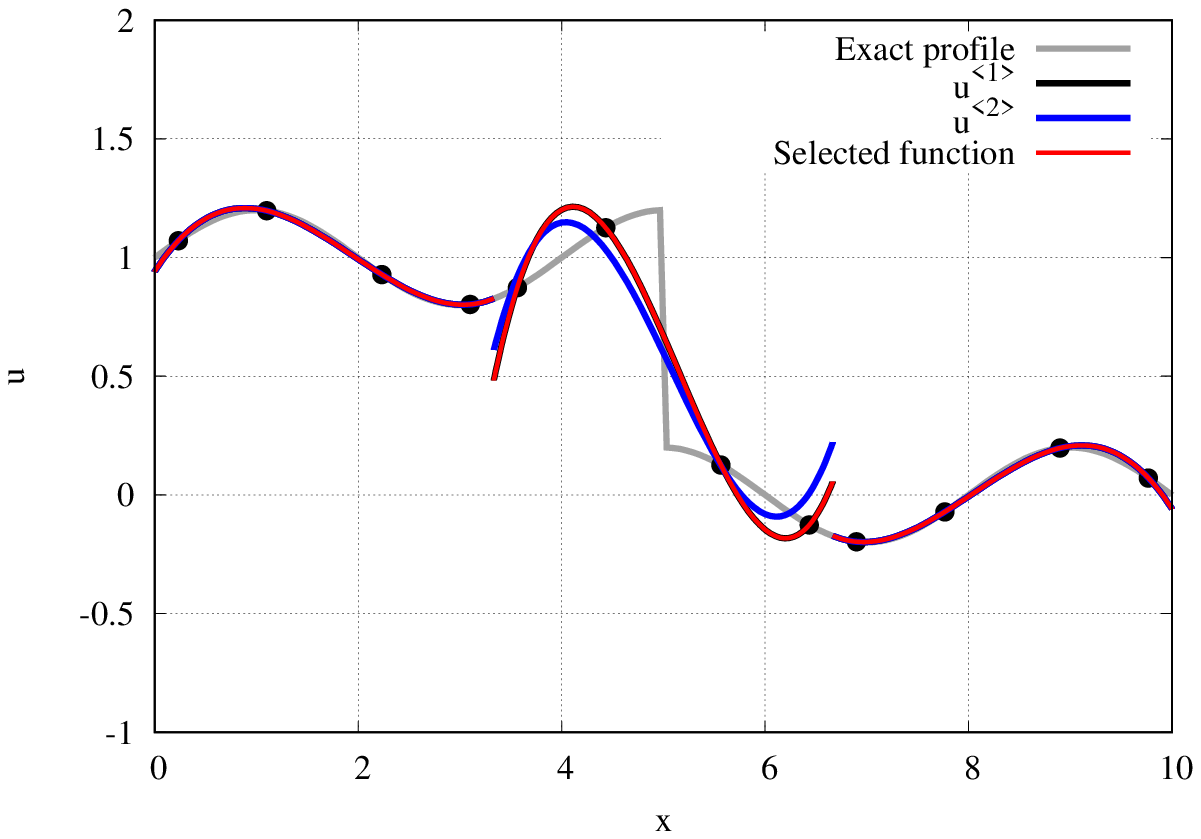}}
\caption{Polynomial selection for the Case 3 solution (smooth and discontinuous flow given by Eq.\eqref{eq_Case3}); polynomial degree $K=3$.
The grey lines show the analytical solution;
black lines show the high-order polynomial approximation $u^{<1>}_i$;
blue lines show the stable function $u^{<2>}_i$;
red lines show the selected function based on TVB or BVD criteria.
The computational domain is divided into three cells, where the functions in the first and third cell are fixed to $u^{<1>}_i$
 }\label{fig_Case3_K3}
\end{figure}

\section{Linear advection equation: $K=3$ cases}
In this section, the cases with $K=3$ in the numerical test for the linear advection equation are presented for the simple WENO reconstruction.
The results are quite encouraging, where the BVD algorithm provides the better results than the TVB limiter cases both with the small and large $M$ values.
\renewcommand{\datahead}{bv_data_cfl0p1_lhs3_rhs1_cor0_sp0_init}
\renewcommand{\timestep}{00001440}
\renewcommand{\wid}{0.3275} 
\begin{figure}[!htbp]
    \centering
    \setcounter{subfigure}{0}
    \subfigure[][{\scriptsize BVD/SWENO}]{\includegraphics[width=\wid\textwidth,clip]{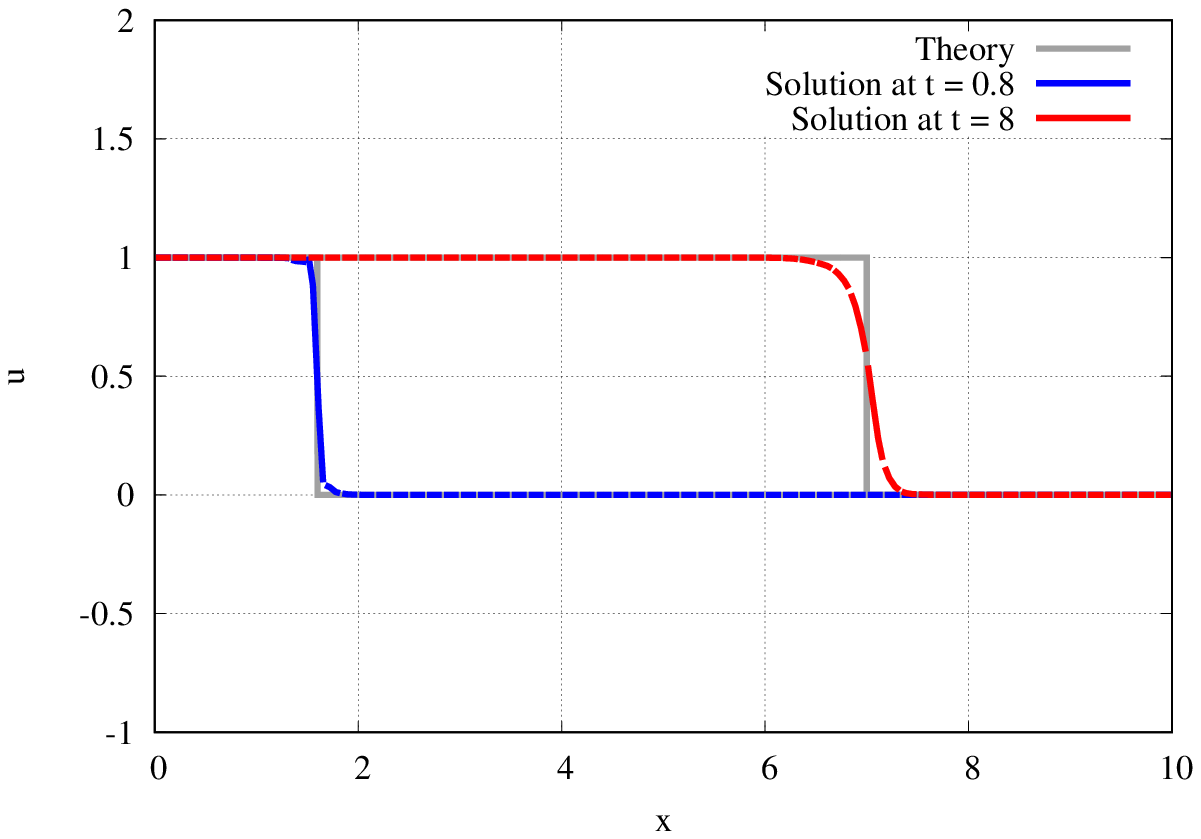}}
    \subfigure[][{\scriptsize TVB/SWENO $(M=0.0)$}]{\includegraphics[width=\wid\textwidth,clip]{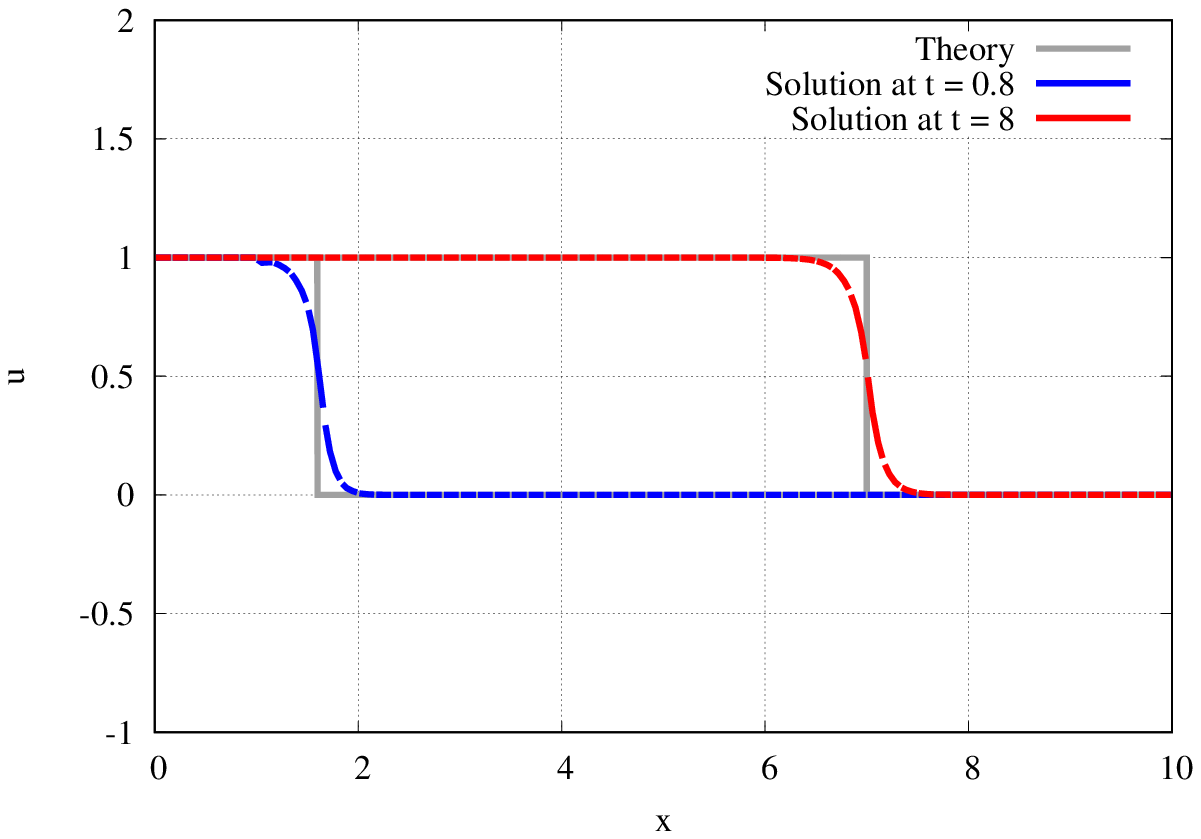}}
    \subfigure[][{\scriptsize TVB/SWENO $(M=200)$}]{\includegraphics[width=\wid\textwidth,clip]{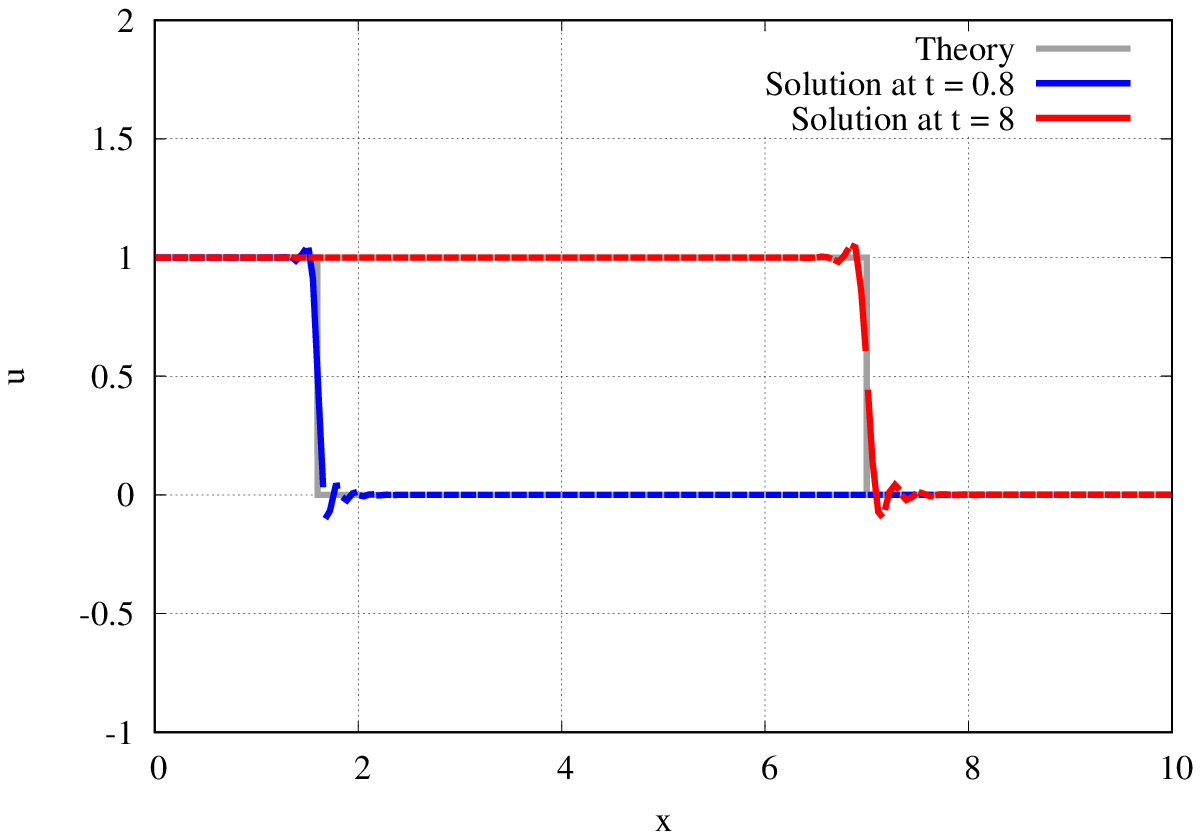}}
\caption{Computational results of a linear advection equation for the case 5 initial condition (discontinuous flow given by Eq.\eqref{eq_Case5}); the total DoF is $240$ with a polynomial degree $K=3$.
The grey lines show the theoritical solution;
blue lines show the solution $u(x,t)$ at $t=0.8$;
red lines show the solution $u(x,t)$ at $t=8.0$.
 }\label{fig_Case5_K3}
\end{figure}
\begin{figure}[!htbp]
    \centering
    \setcounter{subfigure}{0}
    \subfigure[][{\scriptsize BVD/SWENO}]{\includegraphics[width=\wid\textwidth,clip]{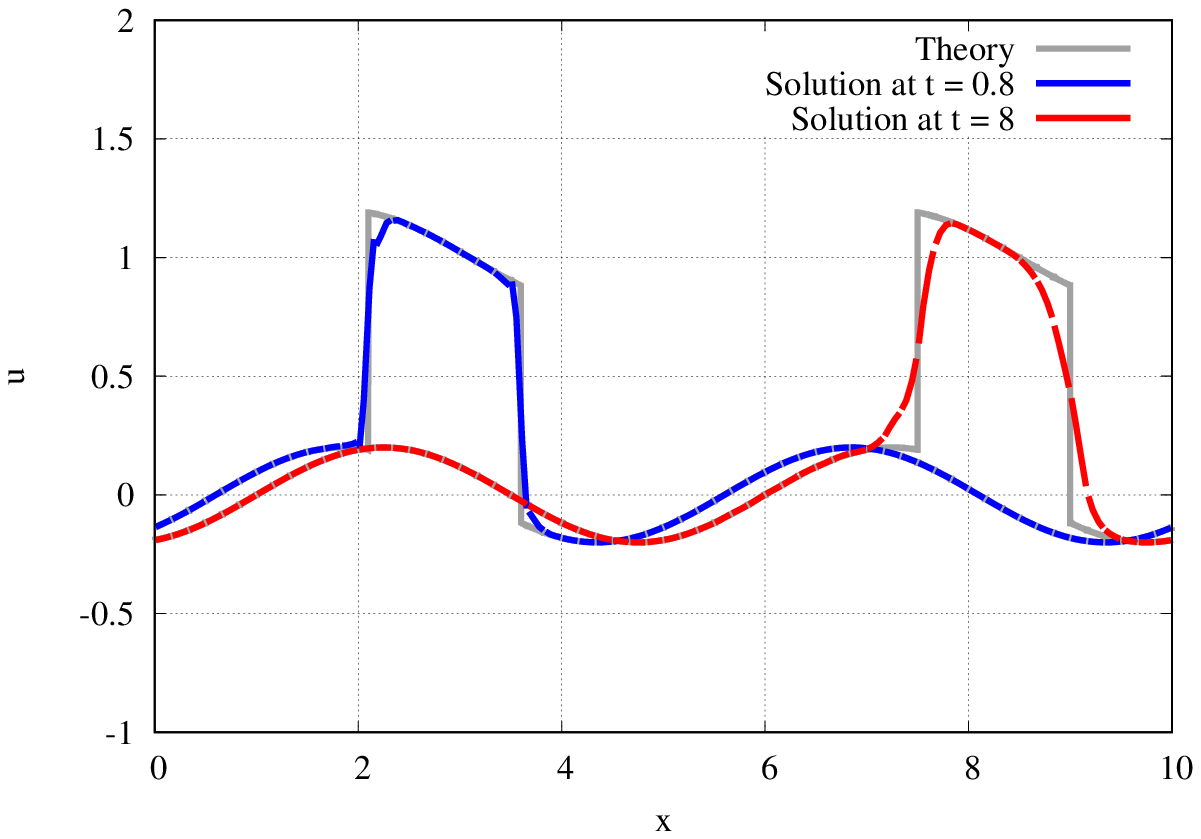}}
    \subfigure[][{\scriptsize TVB/SWENO $(M=0.0)$}]{\includegraphics[width=\wid\textwidth,clip]{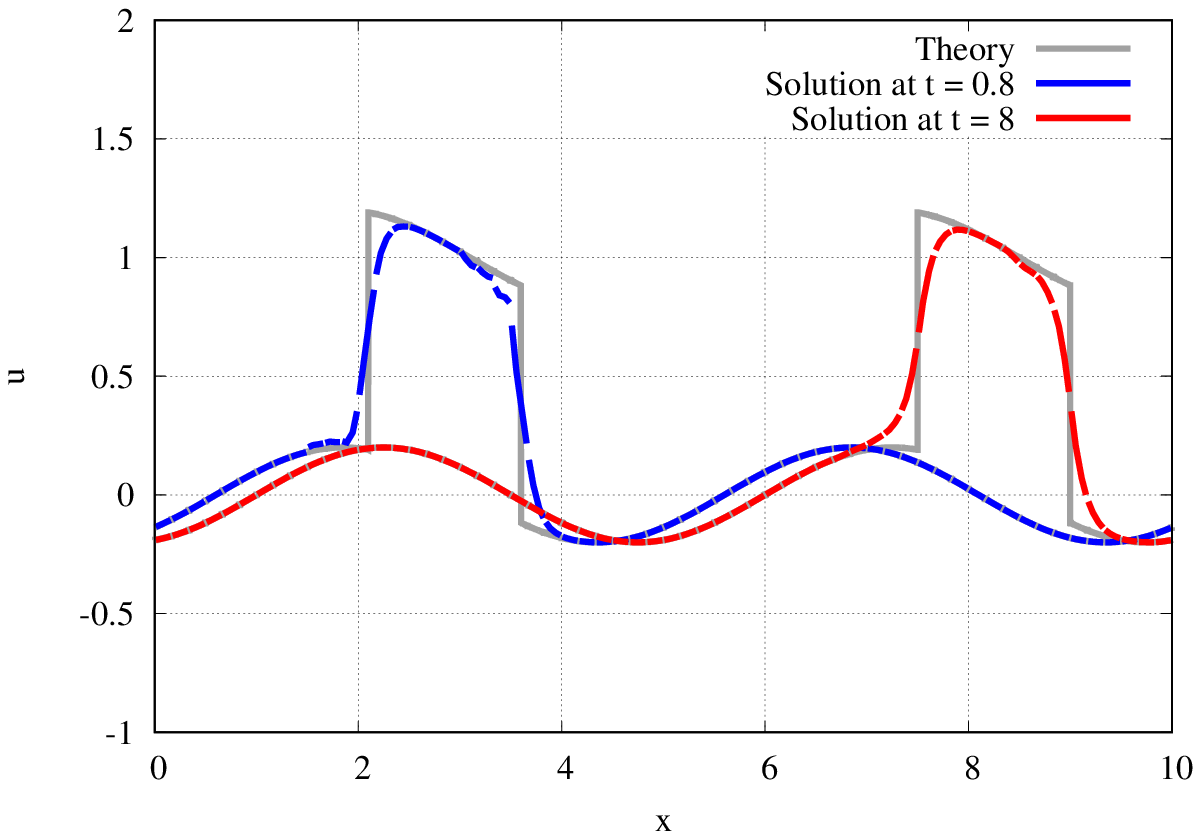}}
    \subfigure[][{\scriptsize TVB/SWENO $(M=200)$}]{\includegraphics[width=\wid\textwidth,clip]{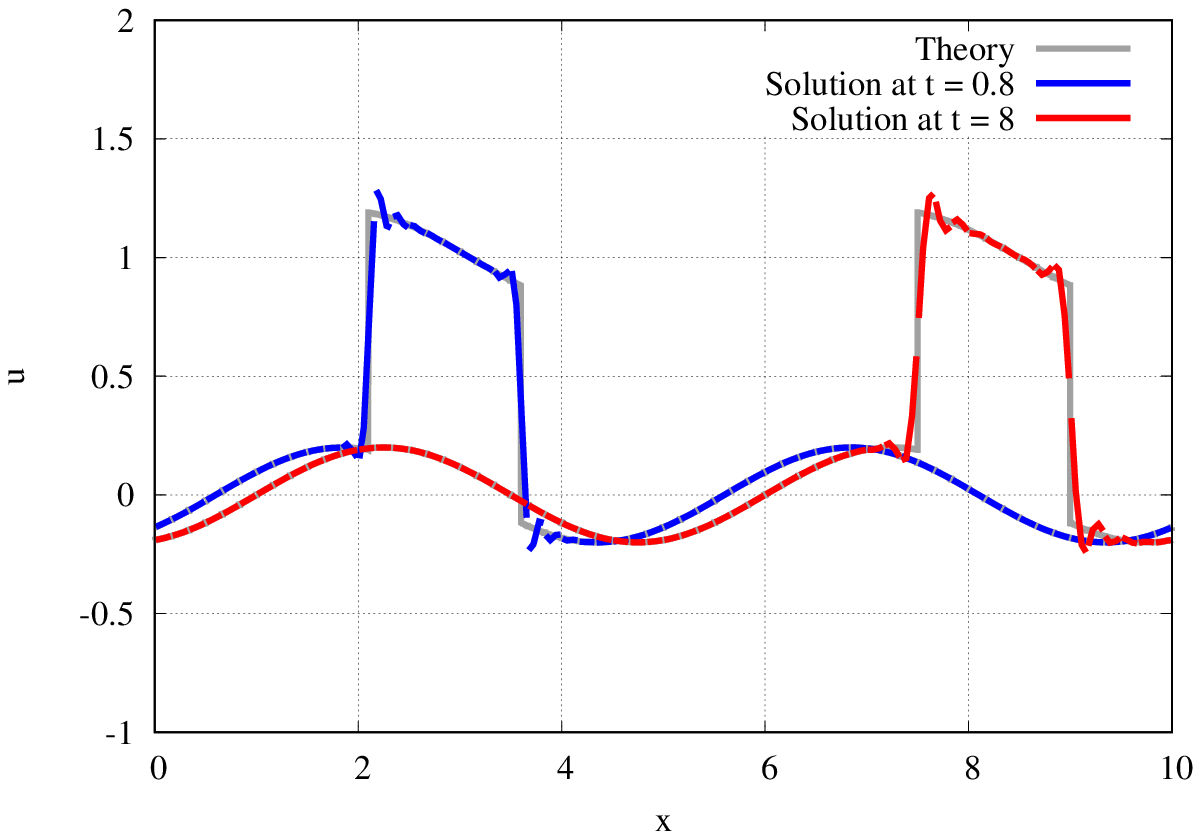}}
\caption{Computational results of a linear advection equation for the case 6 initial condition (convex, smooth, and discontinuous flow given by Eq.\eqref{eq_Case6}); the total DoF is $240$ with a polynomial degree $K=3$.
The grey lines show the theoritical solution;
blue lines show the solution $u(x,t)$ at $t=0.8$;
red lines show the solution $u(x,t)$ at $t=8.0$.
 }\label{fig_Case6_K3}
\end{figure}

\end{document}